\documentclass[reqno,twoside,11pt]{amsart}
\usepackage{amsmath, amsfonts, amssymb, esint, amsthm, nicefrac, bbm, dsfont, enumitem}
\usepackage{epsfig, graphicx}
\usepackage{hyperref, cite, comment}

\def\sideremark#1{\ifvmode\leavevmode\fi\vadjust{\vbox to0pt{\vss
 \hbox to 0pt{\hskip\hsize\hskip1em
\vbox{\hsize2cm\small\raggedright\pretolerance10000
 \noindent #1\hfill}\hss}\vbox to8pt{\vfil}\vss}}}

\setlist[itemize]{leftmargin=*}
\DeclareMathAlphabet{\mathpzc}{OT1}{pzc}{m}{it}
\DeclareMathOperator{\dist}{dist}
\newcommand{\Eh}{\mathcal{E}^h}
\newcommand{\dx}{\;\mathrm{d}x}

\newcommand{\R}{\mathbb{R}}

\newcommand{\G}{g}
\newcommand{\Gc}{\overset{\Gamma}\longrightarrow}
\def\endproof{\hspace*{\fill}\mbox{\ \rule{.1in}{.1in}}\medskip }

\renewcommand{\epsilon}{\varepsilon}

\numberwithin{equation}{section}

\theoremstyle{plain}
\newtheorem{theorem}{Theorem}[section]

\newtheorem{corollary}[theorem]{Corollary}
\newtheorem{proposition}[theorem]{Proposition}

\theoremstyle{definition}
\newtheorem{remark}[theorem]{Remark}

\newtheorem{problem}[theorem]{Open Problem}

\begin{document}
\title[Geometry, Analysis and Morphogenesis]
{Geometry, Analysis and Morphogenesis: Problems and Prospects} 
\author{Marta Lewicka and L. Mahadevan}
\address{M. Lewicka: University of Pittsburgh, Department of Mathematics, 139 University Place, Pittsburgh, PA 15260, USA} 
\address{L. Mahadevan:  School of Engineering and Applied Sciences, and Departments of Physics, and Organismic and Evolutionary Biology, Harvard University, Cambridge, MA 02138}
\email{lewicka@pitt.edu, lmahadev@g.harvard.edu} 

\begin{abstract}
The remarkable range of biological forms in and around us, such as the undulating shape of a leaf or flower in the garden, the coils in our gut, or the folds in our brain, raise a number of questions at the interface of biology, physics and mathematics. How might these shapes be predicted, and how can they eventually be designed? We review our current understanding of this problem, that brings together analysis, geometry and mechanics in the description of the morphogenesis of low-dimensional objects. Starting from the view that shape is the consequence of metric frustration in an ambient space, we examine the links between the classical Nash embedding problem and biological morphogenesis. Then, motivated by a range of experimental observations and numerical computations, we revisit known rigorous results on curvature-driven patterning of thin elastic films, especially the asymptotic behaviors of the solutions as the (scaled) thickness becomes vanishingly small and the local curvature can become large. 
Along the way, we discus open problems that include those in mathematical modeling and analysis along with questions driven by the allure of being able to tame soft surfaces for applications in science and engineering.  
\end{abstract}

\date{\today}
\maketitle

\section{Introduction}

A walk in the garden, a visit to the zoo, or watching a nature documentary reminds us of the remarkable range of living forms on our planet. How these shapes come to be is a question that has interested scientists for eons, and yet it is only over the last century that we have finally begun to grapple with the framework for morphogenesis, a subject that naturally brings together biologists, physicists and mathematicians. This confluence of approaches is the basis for a book, equally lauded for both its substance and its scientific style, D'Arcy Thompson's opus \lq\lq On growth and form" \cite{Th} where the author says: {\em \lq\lq An organism is so comp\-lex a thing, and growth so complex a phenomenon, that for  growth to be so uniform and constant in all the parts as to keep the whole shape unchanged would indeed be an unlikely and an unusual circumstance. Rates vary, proportions change, and the whole configuration alters accordingly."} 

From a mathematical and mechanical perspectives, this suggests a simple principle: differential growth in a body leads to residual strains that will generically result in  changes in the shape of a tissue, organ or body. Surprisingly then, it is only recently that this principle has been taken up seriously by both experimental and theoretical communities as a viable candidate for patterning at the cellular and tissue level, perhaps because of the dual difficulty of measuring, and calculating the mechanical causes and consequences of these effects. Nevertheless, with an increasing number of testable predictions and high throughput imaging in space-time, this geometric and mechanical perspective on morphogenesis has begun to be viewed as a complement to the biochemical aspects of morphogenesis, as famously exemplified by the work of Alan Turing in his prescient paper \lq\lq The chemical basis for morphogenesis" \cite{Turing}. It is worth pointing out that differential diffusion and growth are only parts of an entire spectrum of mechanisms involved in morphogenesis that include differential adhesion, differential mobility, differential affinity and differential activity, all of which we must eventually come to grips with to truly understand the development and evolution of biological shape. 

In this review, we consider the interplay between geometry, analysis and morphogenesis of thin surfaces driven by three motivations:  the allure of quantifying the aesthetic seen in examples such as flowers, the hope of explaining the origin of shape in biological systems, and the promise of mimicking them in artificial  systems  \cite{27, Sharon2}.  While these issues also arise in three-dimensional tissues in such examples as the folding of the brain \cite{brain1,brain2} or the looping of the gut \cite{gut1,gut2}, the separation of scales in slender structures that grow in the plane and out of it links the physical problem of growing elastic films to the geometrical problem of determining a slowly evolving approximately two-dimensional film in three dimensions. Indeed, as we will see, many of the questions we review here are related to a classical theme in differential geometry - that of embedding a shape with a given metric  in a space of possibly different dimension \cite{Nash1, Nash2}, and eventually that of designing  the metric to achieve any given shape. However, the goal now is not only to state the conditions when it might be done (or not), but also to determine the resulting shapes in terms of an appropriate mechanical theory, and understand the limiting behaviors of the solutions as a function of the geometric parameters. 
\begin{figure}
    \centering
    \includegraphics[width=1.1\textwidth]{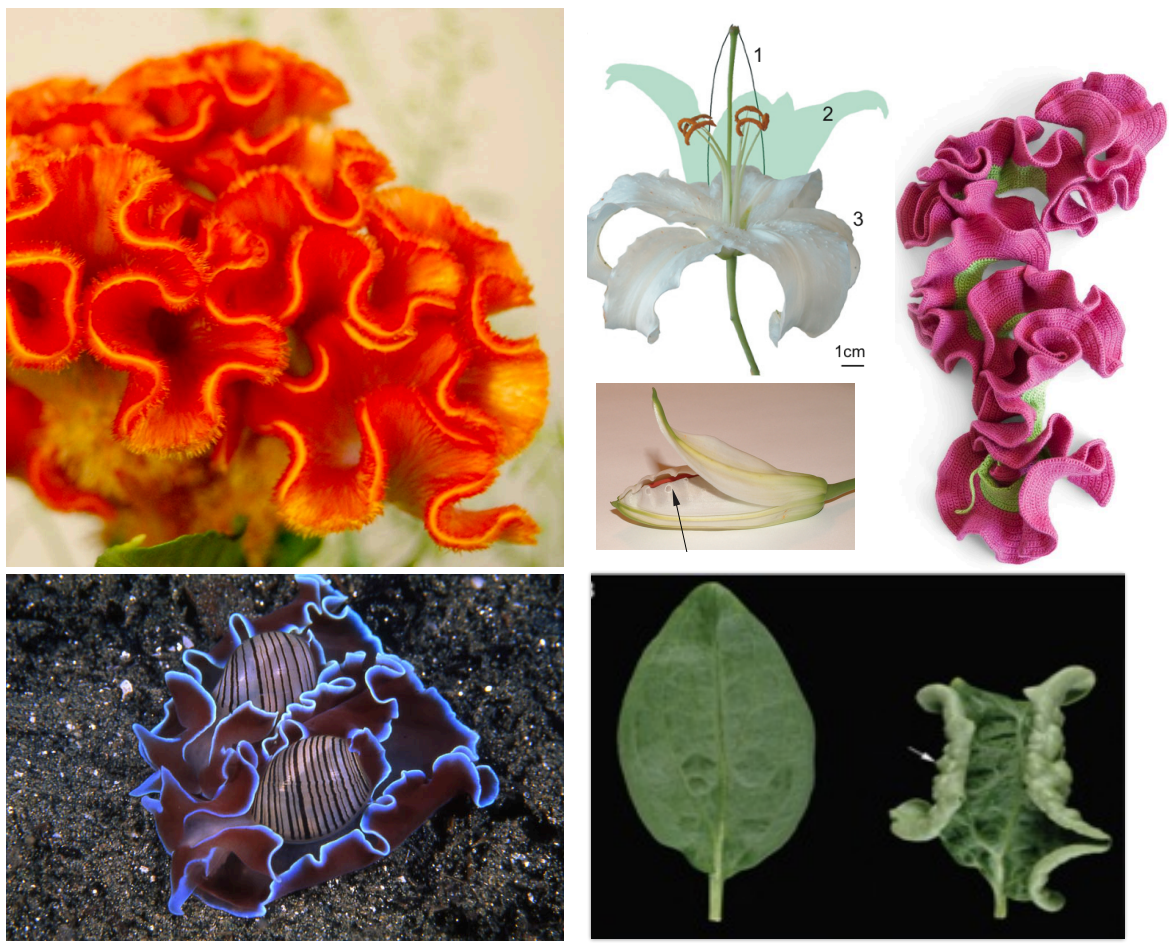}
      \vskip -0in
    \caption{Patterns in a range of systems in terrestrial and aquatic environments show the myriad forms that reflect the consequence of inhomogeneous growth of a thin sheet: the impossibility of embedding an arbitrary biological growth metric coupled with the constraint of minimizing an elastic energy leads to frustration embodied as shape. The examples shown are: (a) a terrestrial cockscomb flower, (b) a marine nudibranch sea-slug, (c) a lily flower in its bud and opened states \cite[Copyright (2021) National Academy of Sciences]{LM2}, (d) a normal and mutant snapdragon leaf, (e) a crocheted scarf. All these are frustrated embeddings of a hyperbolic metric into $\mathbb{R}^3$.}
    \label{fig1}
\end{figure}

\smallskip

The outline of this paper is as follows. Starting from the view that shape is the consequence of metric frustration in an ambient space, in section \ref{sec_new2} we describe the background and objectives of the non-Euclidean elasticity formalism as well as present an example of growth equations in this context. In section \ref{sec_sf}
we examine the links between the classical Nash embedding problem and biological morphogenesis. Then, motivated by a range of experimental observations and numerical computations, we revisit known rigorous results on curvature-driven patterning of thin elastic films in section \ref{sec-3}, where we also offer a new estimate regarding the scaling of the non-Euclidean energies from convex integration. In section \ref{sec4}, we focus on the asymptotic behaviors of the solutions as the (scaled) thickness of the films becomes vanishingly small and the local curvature can become large. In section \ref{sec-5} we digress to consider the weak prestrains and the related Monge-Amp\`ere constrained energies. In section \ref{sec_clas}, the complete range of results is compared with the hierarchy of classical geometrically nonlinear theories for elastic plates and shells without prestrain.  Along the way and particularly in section \ref{sec_fut}, we discuss open mathematical problems and future research directions.

\medskip

\noindent {\bf Acknowledgement.} M. Lewicka was partially supported by NSF grant DMS 2006439. L Mahadevan was partially supported by NSF grants BioMatter DMR 1922321 and MRSEC DMR 2011754 and EFRI 1830901.

\medskip

\section{Non-Euclidean elasticity and an example of growth equations}\label{sec_new2}

An inexpensive surgical experiment serves as a clue to the biological processes at work in determining shape: if one takes a sharp knife and cuts a long, rippled leaf into narrow strips parallel to the midrib, the strips flatten out when "freed" from the constraints of being contiguous with each other. This suggests that the shape is the result of geometric frustration and feedback, driven by the twin effects of: embedding a non-Euclidean metric due to inhomo\-gene\-ous growth, and minimizing an elastic energy that selects the particular observed shape. 
Experiments confirm the generality of this idea in a variety of situations, ranging from undulating submarine avascular algal blades, to saddle-shaped, coiled or edge-rippled leaves of many terrestrial plants \cite{Koehl2008,18b}.
Understanding the origin of the morphologies of slender structures as a consequence of either their growth or the constraints imposed by external forces, requires a mathematical theory for how shape is generated by inhomogeneous growth in a tissue.

\subsection{Non-Euclidean elasticity}
Biological growth arises from changes in four fields: cell number, size, shape and motion, all of which conspire to determine the local metric, which in general will not be compatible with the existence of an isometric immersion. For simplicity, growth is often coarse-grained by averaging over cellular details, thus ignoring microscopic structure due to cellular polarity, orientation (nematic order), anisotropy etc. While recent work has begun to address these more challenging questions \cite{Marchetti-RMP, Warner-ARCMP}, we limit our review to the case of homogeneous, isotropic thin growing bodies. This has proceeded along three parallel paths, all leading to a set of coupled hyperelliptic pd-es that follow from a variational principle: \begin{enumerate}[leftmargin=7mm]
\item[(i)] by using the differential geometry of surfaces as a starting point to determine a plausible class of elastic energies written in terms of the first and second fundamental forms, or their discrete analogs and deviations from some natural state \cite{Grinspun};
\item[(ii)] by drawing on an analogy between growth and thermoelasticity \cite{Mansfield} and plasticity \cite{Lee}, since they both drive changes in the local metric tensor and the second fundamental form, and using this to build an energy functional whose local mini-ma determines shape;
\item[(iii)] by starting from a three-dimensional theory for a growing elastic body with geometrically incompatible growth tensors, driving the changes of the first and second fundamental forms of a two-dimensional surface embedded in three dimensions \cite{Ciarlet, jed1}. 
\end{enumerate}
The resulting shape can be seen as a consequence of the heterogeneous incompatibility of strains that leads to geometric (and energetic) frustration. This coupling between residual strain and shape implies an energetic formulation of {\em non-Euclidean elasticity} that attempts to minimize an appropriate energy associated with the frustration between the induced and intrinsic geometries. Within this framework, a few different types of problems may be posed: 
\begin{enumerate}[leftmargin=7mm]
\item[(i)] questions about 
the nature of the (regular and singular) solutions that arise; \item[(ii)]  questions on their connection to experimental observations; 
\item[(iii)] problems related to the limiting behavior of the solutions and their associated energies in the limit of small (scaled) thickness; 
\item[(iv)] questions on identifying the form of feedback laws linking growth to shape that lead to the self-regulated reproducible forms seen in nature; 
\item[(v)] problems in formulating of inverse problems in the context of shaping sheets for function. 
\end{enumerate}

\subsection{An example of growth equations}
To get a glimpse of the analytical structures to be investigate, we begin by writing down a minimal theory that couples growth to the shape of a thin lamina of uniform thickness \cite{jed1, Mansfield,18b}, now generalized to account for differential growth:
\begin{equation}\label{Goveq0}
\begin{split}
& \Delta ({\rm tr}\,\boldsymbol \sigma) +\frac{\alpha}{2} \det \boldsymbol \kappa=- \alpha \Delta ({\rm tr}\, \mathbf{s})\\
& \beta\Delta ({\rm tr} \,\boldsymbol \kappa) - {\rm tr}\, (\boldsymbol \sigma \boldsymbol \kappa)=-\beta \Delta ({\rm tr}\, \mathbf{b}). 
\end{split}
\end{equation}
Here, $\Delta$ is the two-dimensional Laplace-Bertram operator, $\boldsymbol \sigma$ is the two-dimensional depth averaged stress tensor, and $\boldsymbol \kappa$ is the curvature tensor. The scalar coefficients $\alpha$ and $\beta$ characterize the elastic moduli of the sheet, assumed to be made of a linear isotropic material; $\alpha$ is the resistance to  stretching (and shearing) in the plane, and $\beta$ is the resistance to bending out of the plane. The right hand side of (\ref{Goveq0}) quantifies the source that drives in-plane differential growth due to a prescribed metric tensor $\mathbf{s}$,  and the out-of-plane differential growth gradient across the thickness due to a prescribed second fundamental form (a curvature tensor) $\mathbf{b}$. 

The first equation in the system (\ref{Goveq0}) corresponds to the incompatibility of the in-plane strain due to both the Gauss curvature and the additional contribution from in-plane differential growth, and it is a geometric compatibility relation. The second equation in the system (\ref{Goveq0}) is a manifestation of force balance in the out-of-plane direction due to the in-plane stresses in the curved shell, and the growth curvature tensor associated with transverse gradients that leads to an effective normal pressure. We observe that $\beta/\alpha = \mathcal{O}(h^2)$, where $h$ is the thickness of the tissue, so there is a natural small parameter in the problem $a = h/L \ll 1$, where $L$ is the lateral size of the system. The nonlinear hyperelliptic equations (\ref{Goveq0}) need to be complemented with an appropriate set of boundary conditions on some combination of the displacements, stresses and their derivatives. However, it is not even clear if and when it is possible to realize reasonable physical surfaces for arbitrarily prescribed tensors $\mathbf{s}, \mathbf{b}$, and so one must resort to a range of approximate methods to determine the behavior of the solutions in general. 

There are two large classes of problems associated with the appearance of fine scales or sharp localized conical features, and characterized by two distinguished limits of (\ref{Goveq0}). These correspond to the situation when either the in-plane stress is relatively large or when it is relatively small.
In the first case, when $| \boldsymbol \sigma| \sim \alpha$, i.e. the case where stretching dominates, one can rescale equations (\ref{Goveq0}) so that they yield the singularly perturbed limit:
\begin{equation}\label{Goveq2}
\begin{split}
& \Delta ({\rm tr}\,\boldsymbol \sigma) +\frac{1}{2} \det \boldsymbol \kappa=- \Delta ({\rm tr}\, \mathbf{s})\\
& a^2\Delta ({\rm tr}\, \boldsymbol \kappa)- {\rm tr}\, (\boldsymbol \sigma \boldsymbol \kappa)=-a^2 \Delta ({\rm tr}\, \mathbf{b}).
\end{split}
\end{equation}
As $a^2 \rightarrow 0$, at leading order the second of the equations above implies ${\rm tr}\, (\boldsymbol \sigma \boldsymbol \kappa) =0$, which has a simple geometric interpretation. Namely, the stress-scaled mean curvature vanishes, which is an interesting generalization of the Plateau-Lagrange problem for minimal surfaces. Then, the system (\ref{Goveq2})  describes a finely decorated minimal surface, where wrinkling patterns appear in regions with a sufficiently negative stress.

In the second case,  when the in-plane stress is relatively small $|\boldsymbol \sigma| \sim \alpha a^2$, i.e. the case when bending dominates, one can rescale  (\ref{Goveq0}) to obtain a different singularly perturbed limit:
\begin{equation}\label{Goveq3}
\begin{split}
& a^2\Delta ({\rm tr}\,\boldsymbol \sigma) +\frac{1}{2} \det \boldsymbol \kappa=-   \Delta ({\rm tr}\, \mathbf{s})\\
& \Delta ({\rm tr} \, \boldsymbol \kappa)- {\rm tr} \,(\boldsymbol \sigma \boldsymbol \kappa)=-  \Delta ({\rm tr}\,\mathbf{b}).
\end{split}
\end{equation}
As $a^2 \to 0$, at leading order the first of the equations above yields $ \det \boldsymbol \kappa=-  2 \Delta ({\rm tr}\,\mathbf{s})$, which can be seen as a  Mange-Amp\`ere type equation for the Gauss curvature. Then, the system  (\ref{Goveq3}) describes a spontaneously crumpled, freely growing sheet with conical and ridge-like singularities, similar to the result of many a failed calculation that end up in the recycling bin. 

Adding the growth terms in (\ref{Goveq0}) is however only part of the biological picture, since in general there is likely to be feedback, i.e. just as growth leads to shape, shape (and residual strain) can change the growth patterns. Then, the growth tensors $\mathbf{s}, \mathbf{b}$ must themselves be coupled to the shape of the sheet via additional (dynamical) equations. 

\begin{problem}
The above description follows the one-way coupling of growth to shape and ignores the feedback from the residual strain. It is known that biological mechanisms inhibit cell growth if the cell experiences sufficient external pressure. Although there are proposals for how shape couples back to growth, this remains a largely open question of much current interest in biology, and we will return to this briefly in the concluding sections.
\end{problem}


\section{Shape from geometric frustration in growing laminae}\label{sec_sf}


The variety of forms seen in the three-dimensional shapes of leaves or flowers, reflects their developmental and evolutionary history and the physical processes that shape them, posing many questions at the nexus of biology, physics and mathematics. From a biological perspective, it is known that genetic mutants responsible for differential cell proliferation lead to a range of leaf shapes \cite{Nath-2003,Coen-review}. From a physical perspective, stresses induced by external loads lead to phenotypic plasticity in algal blades that switch between long, narrow blade-like shapes in rapid flow to broader undulating shapes in slow flow \cite{Koehl2008}.  
Similar questions arise from observations of a blooming flower, long an inspiration for art and poetry, but seemingly not so from scientific perspectives.  When a flower blossoms, its petals change curvature on a time scale of a few hours, consistent with the idea that these movements are driven by cellular processes. In flowers that bloom once, differential cell proliferation is the dominant mode of growth, while in those that open and close repeatedly, cell elongation plays an important role. 

\begin{figure}
    \centering
    \includegraphics[width=\textwidth]{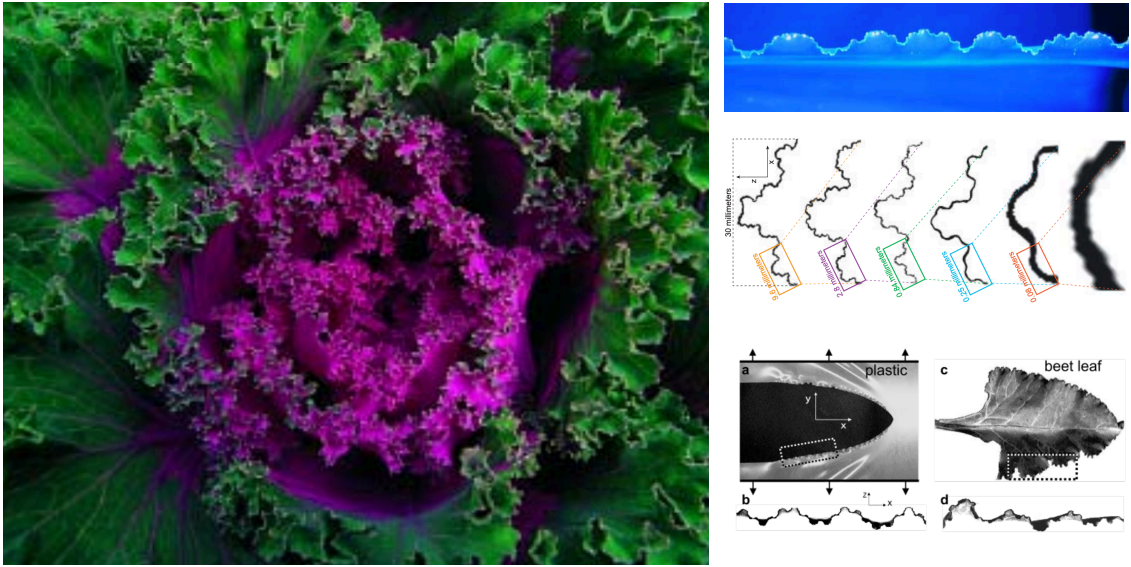}
    \vskip -1.1in
    \caption{Top view of purple kale shows that the edges of the leaf are wrinkled hierarchically, as the thickness of the kale reduces. This can be recapitulated \cite{CRM} in a simple tearing experiment of a garbage bag - the tearing edge thins, is plastic-ally deformed, and thus wrinkles. There is a clear hierarchy seen; further images show the analogy between the mechanically deformed edge of the torn sheet and the edge of a beet leaf.} 
    \label{fig2}
\end{figure}

Although proposed explanations for petal movements posit a difference in growth rate of its two sides (surfaces) or an active role for the midribs, experimental, the theoretical and computational studies \cite{LM2} have shown that the change of the shape of a lamina is due to excess growth of the margins relative to the center (see Figure \ref{fig1}).
Indeed, there is now ample evidence of how relative growth leads to variations in shape in such contexts as leaves, flowers, micro-organisms (i.e. euglenids), swelling sheets of gels, 4d printed structures etc. \cite{26, 22, 22a, 25, 26, 51, 15, DCB, BMT, D}.    A particularly striking example is that of the formation of self-similar wrinkled structures   as shown in the example of a kale leaf in Figure \ref{fig2}. A demonstration of the same phenomenon with everyday materials is also shown in Figure \ref{fig2} - when a garbage bag is torn, its edge shows multiple generations of wrinkles \cite{Sharon2}. 

\subsection{The set-up}

The experimental observations described above suggest a common mathematical framework to understand the origin of shape: an elastic three-dimensional body $\Omega$ seeks to realize a configuration with a prescribed Riemann metric $g$ by means of an isometric immersion. The deviation from or inability to reach such a state, is due to a combination of geometric incompatibility and the requirements of elastic energy minimization. Borrowing from the theory of finite plasticity \cite{Lee}, where a multiplicative decomposition of the deformation gradient into an elastic and a plastic use was postulated, a similar hypothesis was suggested for growth \cite{RHM}, with the underlying hypothesis of the presence of a reference configuration $\Omega$ with respect to which all displacements are measured. 

Let $g:\Omega\to\mathbb{R}^{3\times 3}_{sym, >}$ be a smooth Romanian metric, given on an open, bounded domain $\Omega\subset\mathbb{R}^3$, and let $u:\Omega\to \R^3$ be an immersion that corresponds to the elastic body.    
Excluding nonphysical deformations that change the orientation in any neighborhood of the immersion, a natural way to pose the question of the origin of shape is by postulating that it arises from a variational principle that minimizes an elastic energy $\mathcal{E}(u)$ which measures how far a given $u$ is from being an orientation-preserving realization of $\G$. Equivalently, $\mathcal{E}(u)$ quantifies the total point-wise deviation of $\nabla u$ from $\G^{1/2}$, modulo orientation-preserving rotations that do not cost any energy. The infamy of $\mathcal{E}$ in absence of any forces or boundary conditions is then indeed strictly positive for a non-Euclidean $\G$, pointing to existence of residual strain.

Since the matrix $g(x)$ is symmetric and positive definite, it possesses a unique symmetric, positive definite square root $A(x)={g(x)}^{1/2}\in\mathbb{R}^{3\times 3}_{sym, >}$ which corresponds to the growth pres-train. This allows us to define an energy: 
\begin{equation}\label{nonE}
\mathcal{E}(u) = \int_{\Omega} W\left((\nabla u)A^{-1}\right)~\mbox{d}x \qquad
\forall u\in W^{1,2}(\Omega,\mathbb{R}^3),
\end{equation}
where the energy density $W:\R^{3\times 3}\to [0,\infty]$
obeys the principles of material frame invariance (with respect to the special orthogonal group of proper rotations $SO(3)$), normalization,
non-degeneracy, and material consistency  valid for all
$F\in\mathbb{R}^{3\times 3}$ and all $R\in SO(3)$:
\begin{equation}\label{elastic_dens}
\begin{split}
& W(RF) = W(F),  \qquad 
 W(Id_3) = 0, \qquad W(F) \geq c~\mbox{dist}^2(F, SO(3)), \\
& W(F)\to +\infty \quad \mbox{as } \det F\to 0+, \quad \mbox{ and}\quad \forall \det
F\leq 0 \quad W(F) = +\infty.
\end{split}
\end{equation}
These models\footnote{Examples of  $W$ satisfying these conditions are: 
$W_1(F) = |(F^TF)^{1/2} - \mbox{Id}|^2 + |\log \det F|^q$,  or $W_2(F)=|(F^TF)^{1/2} - \mbox{Id}|^2 + \left|({\det F})^{-1} - 1\right|^q$ for $\det F>0$, where $q>1$ and $W_{1,2}$ equal $+\infty$ if $\det F\leq 0$.}, corresponding to a range of hyperelastic energy functionals that approximate the behavior of a large class of elastomeric materials, are consistent with microscopic derivations based on statistical mechanics, and naturally reduce to classical linear elasticity when $|(F^TF)^{1/2} - \mbox{Id}| \ll 1$. Minimizing the energy (\ref{nonE}) is thus a prescription for shape, and may be defined naturally in terms of the energetic cost of deviating from an isometric immersion.

\subsection{Isometric immersions and residual stress}\label{ss_2.2}

The model in (\ref{nonE}) assumes that the $3$d elastic body $\Omega$ seeks to realize a configuration with 
a prescribed Riemannian metric $g$, through minimizing the
energy that is determined by the elastic part $F_e=(\nabla u)A^{-1}$ of the
deformation gradient $\nabla u$. 
Observe that $W(F_e)=0$ if and only if $F_e\in SO(3)$ in $\omega$, or equivalently when:
$$(\nabla u)^T\nabla u = g\quad\mbox{ and } \quad\det\nabla u > 0 \quad\mbox{in }\;\omega,$$
Further, any $u\in W^{1,2}(\Omega,\R^3)$ that satisfies
the above, must automatically be smooth. Indeed, writing $\nabla u =
Rg^{1/2}$ for some rotation field $R:\Omega\to SO(3)$, it follows that
$u\in W^{1,\infty}$ and so $\mbox{div}\big(\mbox{cof}\,\nabla
u\big)=0$ holds, in the sense of
distributions\footnote{The divergence of a matrix field is taken row-wise.}.
Further, we have:
\begin{equation}\label{iso}
\mbox{cof}\,\nabla u =\mbox{cof}(Rg^{1/2})=\det(Rg^{1/2})
\big(Rg^{1/2}\big)^{T, -1} = \sqrt{\det g}\big(R g^{-1/2}\big)=
\sqrt{\det g}\big((\nabla u) g^{-1}\big).
\end{equation}
It follows that each of the three scalar components of $u$ is harmonic with respect
to the Laplace-Beltrami operator $\Delta_g$ and thus $u$ must be smooth:
$$0=\mbox{div}\big(\mbox{cof}\,\nabla u\big)=\mbox{div}\big((\sqrt{\det g})(\nabla u) g^{-1}\big) =
  \sqrt{\det g}\cdot \big[\Delta_gu^m\big]_{m=1}^3.$$

Thus, $\mathcal{E}(u)=0$ if and only if
the deformation $u$ is an orientation preserving isometric immersion of $g$ 
into $\mathbb{R}^3$. Such smooth (local) immersion exists \cite[Vol. II, Chapter 4]{Spivak}
and is automatically unique up to rigid motions of $\R^3$, if and only if the Riemann curvature tensor $[\mathcal{R}_{ij, kl}]_{i,j,k,l=1\ldots 3}$ of $g$ vanishes identically throughout $\Omega$. 

\begin{figure}[htbp]
\centering
    \includegraphics[width=\textwidth]{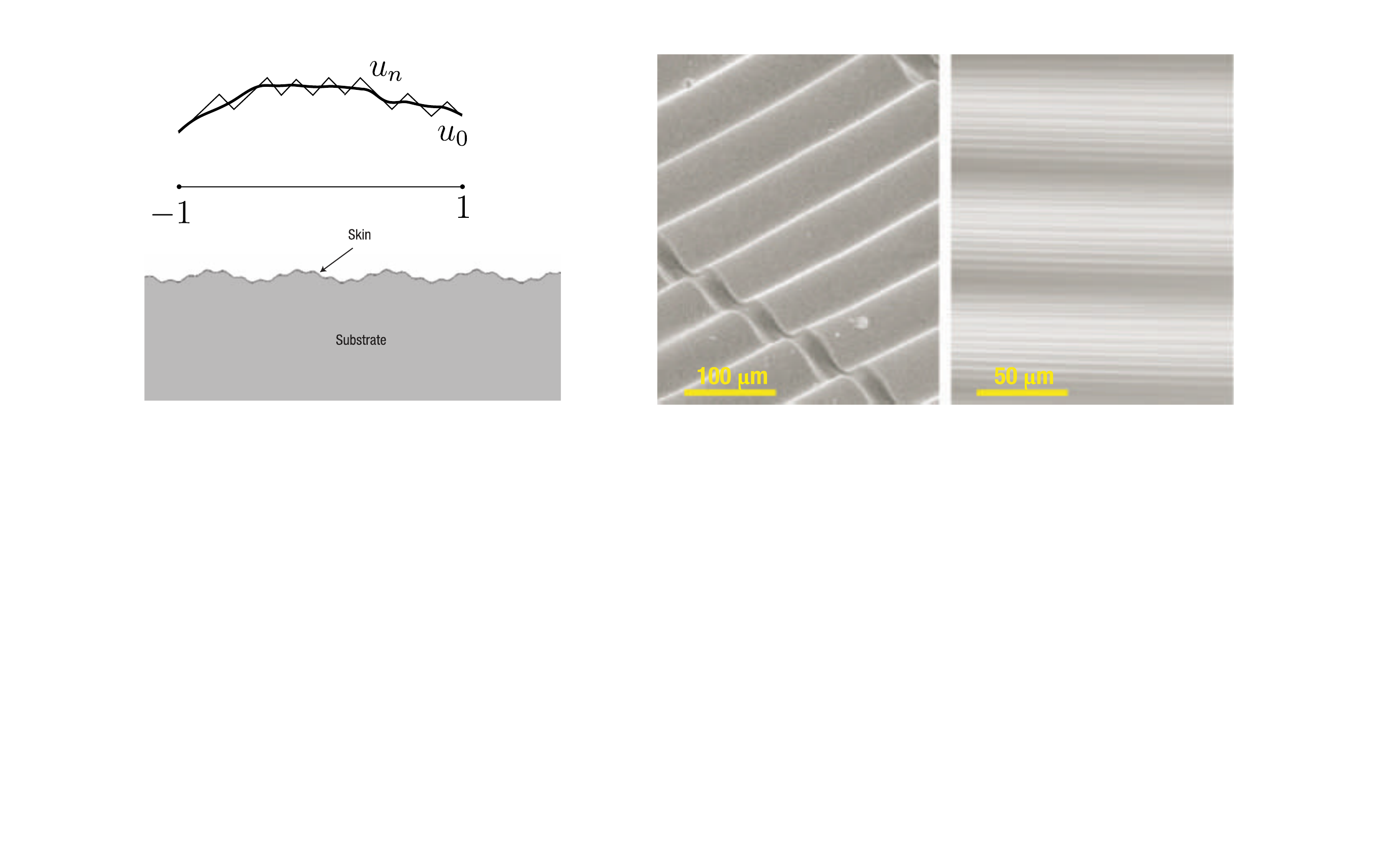}
    \vskip -2in
\caption{(a) A short map approximation of $u_0$ (darker line) by a zigzag $u_n$ with $u_n'=\pm 1$. (b) A computational realization of hierarchical wrinkles that arise when a thin stiff film is coated atop a soft substrate and the system is then subject to a reduction in temperature that leads to differential shrinkage \cite{ERMVMG}. (c) An experimental realization of the hierarchical wrinkles that shows two (of a total of six) generations of wrinkles. The three examples serve to link convex integration to models and experiments in materials science \cite{ERMVMG}.}
\label{zigzag_fig}
\end{figure}

It is instructive to point out that one could define the energy as the difference between the prescribed metric $g$, and the pull-back metric of $u$ on $\Omega$:
\begin{equation*}
I(u)= \int_{\Omega} |(\nabla u)^T\nabla u - g|^2 ~\mbox{d}x.
\end{equation*}
From a variational point of view, the formulation above does not capture an essential aspect of the physics, namely that thin laminae resist bending deformations that are a consequence of the extrinsic geometry, and thus depend on the mean curvature as well. Indeed, the functional $I$ always minimizes to $0$ because there always exists a Lipschitz isometric immersion $u\in W^{1,\infty}(\Omega,\mathbb{R}^3)$ of $g$, for which $I(u) = 0$. If $\mathcal{R}_{ij,kl}(x)\neq 0$ for some $x\in\Omega$, then such $u$ must have a {\em folding structure} \cite{gromov} around $x$; it cannot be orientation preserving (or reversing) in any open neighborhood of $x$.
Perhaps even more surprisingly, the set of such Lipschitz isometric immersions is dense in the set of short immersions: for every $u_0\in \mathcal{C}^1(\bar\Omega,
\R^3)$ satisfying $(\nabla u_0)^T\nabla u_0<g$,\footnote{That is, the matrix $g(x)
  - \nabla u_0(x)^T\nabla u_0(x)$ is positive definite at each $x\in\omega$.} there exists a sequence
$\{u_n\in W^{1,\infty}(\Omega, \R^3)\}_{n\to\infty}$ satisfying:
$$I(u_n) = 0 \quad\mbox{ and } \quad \|u_n-u_0\|_{\mathcal{C}(\Omega)}
\to 0\quad \mbox{ as }\; n \to\infty.$$
The above statement is an example of the {\em $h$-principle} in differential geometry and it follows through the method of convex integration (the Nash-Kuiper scheme), to which we come back in the following sections. An intuitive example in dimension $1$ is shown in Figure \ref{zigzag_fig}.
Setting $g\equiv 1$ on $\Omega=(-1,1)\subset\R^1$, it is easily seen that any $u_0:(-1,1)\to\R$ with Lipschitz constant less than $1$ can be uniformly approximated by $u_n$ having the form of a zigzag, where $|u_n'|=1$.

Regarding the energy $\mathcal{E}$ in (\ref{nonE}), in \cite{non-Euclid} it has been proved that $\inf\mathcal{E}>0$ for any $g$ with no orientation-preserving isometric immersion. This results in the dichotomy: either $g$ and $\mathcal{E}$ are, by a smooth change of variable, equivalent to the case with $g=Id_3$ and $\min \mathcal{E}=0$, or otherwise the zero energy level cannot be achieved through a sequence of weakly regular $W^{1,2}$ deformations. The latter case points to existence of {\em residual strain at free equilibria}. 


\begin{proposition}\label{basic}\cite{non-Euclid}
If $[\mathcal{R}_{ij,kl}] \not\equiv 0$ in $\Omega$, then
$\inf\big\{\mathcal{E(}u); ~u\in W^{1,2}(\Omega,\mathbb{R}^3)\big\}>0$.
\end{proposition}
\noindent {\it Sketch of proof.}
Assume, by contradiction, that $\mathcal{E}(u_n)\to 0$ along some
sequence $\{u_n\in W^{1,2}(\Omega, \R^3)\}_{n\to\infty}$. By
truncation and approximation in Sobolev spaces, we may without
loss of generality assume that each $u_n$ is Lipschitz with a
uniform Lipschitz constant $M$. Decompose $u_n=z_n+w_n$ as a sum of a deformation that is clamped at the boundary:
$$[\Delta_g z_n^m]_{m=1}^3 = [\Delta_g u_n^m]_{m=1}^3 \quad \mbox{in }\Omega, \qquad z_n=0 \quad\mbox{on } \partial\Omega, $$
and a harmonic correction: $[\Delta_g w_n^m]_{m=1}^3 = 0$ in $\Omega$,
with $w_n=u_n$ on $\partial\Omega$. Observe that:
\begin{equation*}
\begin{split}
\int_\Omega \big\langle (\sqrt{\det g})&\nabla z_n g^{-1}:\nabla z_n\big\rangle\dx
= \int_\Omega \big\langle (\sqrt{\det g})\nabla u_n g^{-1}:\nabla z_n\big\rangle\dx
\\ & = \int_\Omega \big\langle (\sqrt{\det g})\nabla u_n g^{-1}:\nabla
z_n\big\rangle - \langle \mbox{cof}\,\nabla u_n : \nabla z_n\rangle \dx \\ & \leq
\|\nabla z_n\|_{L^2} \big\| (\sqrt{\det g})\nabla u_n g^{-1} -\mbox{cof}\,\nabla u_n\big\|_{L^2},
\end{split}
\end{equation*}
where the first equality follows by $z_n=0$ on
$\partial\Omega$, as $\mbox{div}\big((\sqrt{\det g})(\nabla u) g^{-1}\big) = \sqrt{\det g} \big[\Delta_gu^m\big]_{m=1}^3$, while the second by $\mbox{div}(\mbox{cof}\,\nabla u)=0$. The left hand side is also equivalent to $\|\nabla z_n\|_{L^2}^2$, so:
\begin{equation}\label{ma1}
\|\nabla z_n\|_{L^2} \leq C\big\| (\sqrt{\det g})\nabla u_n g^{-1} -\mbox{cof}\,\nabla u_n\big\|_{L^2}
\leq C_M \mathcal{E}(u_n)^{1/2}\to 0.
\end{equation}
Above, we used (\ref{iso}) that assures vanishing of the expression under the norm when $(\nabla u_n)A^{-1}\in SO(3)$, together with Lipschitz continuity of the operator in the integral expression for $\mathcal{E}$. In particular, we get that both sequences $\{\nabla z_n\}_{n\to \infty}$ and $\{\nabla w_n\}_{n\to\infty}$ are bounded in $L^2$.

Since $\{w_n\}$ are harmonic, this further implies that $\nabla u_n$ converges, up to a subsequence, in $L^2_{loc}$ to some $\nabla u$. By (\ref{ma1}) then: $\nabla u_n\to\nabla u$, which yields $\mathcal{E}(u) =0$ and ends the proof.
\endproof

\begin{problem}
In the above context, prove that $\inf \mathcal{E}$ as in Proposition \ref{basic} is equivalent to $\|[\mathcal{R}_{ij, kl}]\|_{H^{-2}}^{1/2}$, up to multiplicative constants depending on $\Omega$ and $W$ but not on $g$. The case of $\Omega\subset\mathbb{R}^2$ and $\mathcal{R}$ replaced by the Gaussian curvature has been considered in \cite{KS}.
\end{problem}

\section{Microstructural patterning of thin elastic prestrained films}\label{sec-3}



\begin{figure}[htbp]
\centering
    \includegraphics[width=1.02\textwidth]{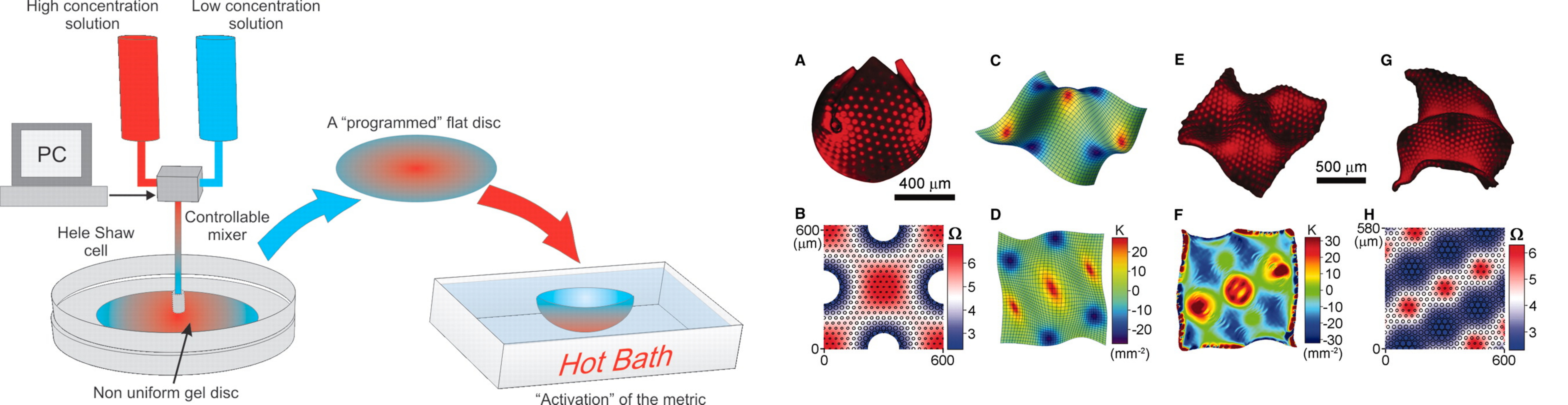}
\caption{Imposing nontrivial target metrics in sheets of NIPA gels. The figures shown are: (a) radially symmetric discs cast by injecting the solution into the gap between two flat glass plates through a central hole  \cite[Reprinted with permission from AAAS]{27}, (b) nonaxisymmetric swelling patterns constructed by half-tone gel lithography in 
\cite[Reprinted with permission from AAAS]{26}.} 
\label{fig5.5}
\end{figure}

Inspired partly by biological observations of growth-induced patterning in thin sheets, and the promise of engineering applications, various techniques have been invented for the construction of self-actuating elastic sheets with prescribed target metrics. The materials typically involve the use of gels that respond to pH, humidity, temperature and other stimuli \cite{Kuma} and result in the formation of complex controllable shapes (see Figure \ref{fig5.5}) that include both large-scale buckling and small-scale wrinkling forms. 

In one example \cite{27},  NIPA monomers with BIS cross-linker in water and a catalyst, leads to the polymerization of a cross-linked elastic hydrogel, which undergoes a sharp, reversible, volume reduction transition at a threshold temperature, allowing for temperature controlled swelling in thin composite sheets.
Another method \cite{26} involves the photopatterning of polymer films to yield temperature-responsive gel sheets with the ability to print nearly continuous patterns of swelling. A third method \cite{22} uses 3d printing of complex-fluid based inks to create bilayers with varying line density and anisotropy, to achieve control over the extent and orientation of swelling.  
All these methods have been used to fabricate surfaces with constant Gaussian curvature (spherical caps, saddles, cones) or zero mean curvature (Enneper’s surfaces), as well as more complex and nearly closed shapes. A natural question that these controlled experiments raise is the ability (or lack thereof) of the resulting patterns to approximate isometric immersions of prescribed metrics. From a mathematical perspective, this leads to questions of the asymptotic behavior of energy minimizing deformations and their associated energetics.

\subsection{The set-up}
In this and the next sections, we will consider a family $(\Omega^h, u^h, \G, A, \Eh)_{h>0}$ (or more generally  $(\Omega^h, u^h, \G^h, A^h, \Eh)_{h>0}$) 
given in function of the film's thickness parameter $h$. The main objective is now to predict the scaling of $\inf \Eh$ as $h\to 0$ and to analyze the asymptotic behavior of minimizing deformations $u^h$ in relation to the curvatures associated with the prestrain $A$.
We assume that $A= g^{1/2}:\overline{\Omega^1}\to\mathbb{R}^{3\times
  3}_{sym, >}$ is a smooth, symmetric and
positive definite tensor field on the unit thickness domain
$\Omega^1$, where for each $h>0$ we define:
\begin{equation*}
\Omega^h=\omega \times \big(-\frac{h}{2}, \frac{h}{2}\big).
\end{equation*}
The open, bounded set $\omega\subset\mathbb{R}^2$
with Lipschitz boundary is viewed as the midplate of the thin film $\Omega^h$, on which we pose the energy of elastic deformations:
\begin{equation}\label{Ih}
\Eh(u^h) = \frac{1}{h}\int_{\Omega^h} W\big((\nabla
u^h)A^{-1}\big)\dx\qquad \mbox{ for all }\; u^h\in W^{1,2}(\Omega^h,\R^3).
\end{equation}

\subsection{Isometric immersions and energy scaling} 

As in section \ref{ss_2.2}, there is a connection between $\inf\mathcal{E}^h$ and existence of isometric immersions, although this case is
 a bit more subtle. In the context of dimension reduction, this connection relays on the isometric
immersions of the midplate metric $g(\cdot,0)_{2\times 2}$  on $\omega$ into $\R^3$, namely parametrised surfaces $y:\omega\to\R^3$ with:
\begin{equation}\label{iso2} (\nabla y)^T\nabla y = g(\cdot,0)_{2\times 2}\quad  \mbox{in } \omega
\end{equation}
It turns out that existence of $y$ with regularity $W^{2,2}$ is equivalent to the vanishing of $\inf \Eh$ of order square in the
film's thickness $h\to 0$. The following result was proved first for $g=g(x')$ in \cite{non-Euclid} and 
then in the abstract setting of Riemannian manifolds in \cite{ KS14}: 

\begin{theorem}\cite{BLS}\label{beta2}
Let $u^h\in W^{1,2}( \Omega^h, \mathbb R^3)$ satisfy  $\Eh(u^h) \le Ch^2$.  Then we have:
\begin{enumerate}[leftmargin=7mm]
\item[(i)] (Compactness). 
There exist $c^h\in \mathbb R^3$ and $R^h\in SO(3)$ such that the rescaled deformations
$y^h(x', x_3) = R^h u^h(x', hx_3)  - c^h$ converge, up to a subsequence in $W^{1,2}(\Omega^1, \R^3)$,  to some $y \in W^{2,2}(\Omega^1, \mathbb R^3)$
depending only on the tangential variable $x'$ and satisfying (\ref{iso2}).

\item[(ii)] (Liminf inequality). There holds the lower bound: 
\begin{equation}\label{I2}
\liminf_{h\to 0} \frac{1}{h^2}\Eh(u^h) \ge {\mathcal I}_{2,g}(y)
 = \frac {1}{24} \int_{\omega} {\mathcal Q}_2 \big(x', (\nabla y)^T
 \nabla \vec b - \frac{1}{2} \partial_3 g(\cdot,0)_{2\times 2}\big)  \, {\rm d}x', 
\end{equation}
where $\mathcal{Q}_2(x',\cdot)$ are nonnegative quadratic
forms given in terms of $D^2W(Id_3)$ (see (\ref{Qform}), and where $\vec b$ satisfies: $\big[\partial_1 y, \partial_2y, \vec
b\big]\in SO(3)g(\cdot,0)^{1/2}$. Equivalently, $\vec b$ is the Cosserat vector comprising the sheer, in addition to the direction $\vec N$ that is normal to the  surface $y(\omega)$:
\begin{equation}\label{Cosserrat_vector} 
\vec b = (\nabla y) g_{2\times 2} ^{-1} \left [ \begin{array}{l} g_{13} \\ g_{23} \end{array} \right ] 
+ \frac{\sqrt{\det g}}{\sqrt{\det g_{2\times 2}}} \vec N, 
\quad \mbox{with: } \,\, \vec N= \frac{\partial_1 y \times \partial_2 y}{|\partial_1 y \times \partial_2 y|}.  
\end{equation}
\end{enumerate}
Moreover, there holds:
\begin{enumerate}[leftmargin=7mm]
\item[(iii)] (Limsup inequality). For all $y\in W^{2,2} (\omega,
  {\mathbb R}^3)$ satisfying (\ref{iso2}) there exists a sequence
  $\{u^h \in W^{1,2}(\Omega^h, \mathbb R^3)\}_{h\to 0}$ for which convergence as in (i) above holds with $c_h=0$, $R^h = {Id}_3$, and:
$$ \lim_{h\to 0} \frac{1}{h^2}\Eh(u^h) = {\mathcal I}_{2,g}(y). $$ 
\end{enumerate}
\end{theorem}

The energy density in (\ref{I2}) is given in terms of a family of quadratic forms $\mathcal{Q}_2(x', \cdot)$, that carry the two-dimensional reduction of the lowest-order nonzero term in the Taylor expansion of $W$ close to its energy well $SO(3)$, namely\footnote{Both $\mathcal{Q}_3$ and all $\mathcal{Q}_2(x',\cdot)$ are nonnegative definite and depend only on the  symmetric parts of their arguments, in view of assumptions on $W$.}:
\begin{equation}\label{Qform}
\begin{split}
& \mathcal{Q}_2(x', F_{2\times 2}) = \min\big\{
  \mathcal{Q}_3\big(g(x',0)^{-1/2}\tilde F g(x',0)^{-1/2}\big);
  ~ \tilde F\in\mathbb{R}^{3\times 3} \mbox{ with }\tilde F_{2\times 2} = F_{2\times 2}\big\},\\
& \mbox{where:} \quad \mathcal{Q}_3(F) = D^2 W(Id_3)(F,F).
\end{split}
\end{equation}

\medskip

From Theorem \ref{beta2}, one can deduce a counterpart of Proposition \ref{basic},
stating an equivalent condition for existence of a $W^{2,2}$ isometric immersion
of a 2-dimensional metric in $\R^3$.

\begin{corollary}\label{isoim}
A smooth metric $\bar g$ on $\bar\omega\subset\mathbb{R}^2$ has an isometric immersion
$y\in W^{2,2}(\omega,\mathbb{R}^3)$ if and only if $\inf \Eh\leq
Ch^2$ for some (equivalently, for any) metric $g$ on $\Omega^1$
with $g(\cdot, 0)_{2\times 2} = \bar g$.
\end{corollary}

The question of existence of local isometric immersions of a given two-dimensional Riemannian manifold into $\mathbb{R}^3$ is a longstanding problem in differential
geometry, its main feature consisting of finding the optimal regularity. By a classical result in \cite{kuiper}, a $\mathcal{C}^1$ isometric embedding can be obtained by means of convex integration. This statement has been improved in \cite{borisov}  to $\mathcal{C}^{1,\alpha}$ regularity for all $\alpha<1/7$ and analytic metrics $\bar g$, in \cite{CDS} to $\mathcal{C}^2$ metrics, and in \cite{DIS} for all $\alpha<1/5$.\footnote{Of interest is also the result in \cite{DI1/2}, stating that for $\alpha>1/2$ the Levi-Civita connection of any isometric immersion is induced by the Euclidean connection, whereas for any $\alpha<1/2$ this property fails to hold.}
This regularity is far from $W^{2,2}$, where information about the second derivatives is also available.
On the other hand, a smooth isometric immersion exists for some special cases, e.g. for smooth $\bar g$ with uniformly positive or negative Gaussian curvatures $\kappa$ on bounded domains in ${\mathbb R^2}$ \cite[Theorems 9.0.1 and 10.0.2]{HH}. Counterexamples to such theories are largely unexplored. By \cite{iaia}, there exists an analytic metric $\bar g$ with nonnegative $\kappa$ on $2$d sphere,
with no local $\mathcal{C}^3$ isometric embedding. However, such metric always admits a $\mathcal{C}^{1,1}$ embedding \cite{GL, HZ}; for a related example see also \cite{pogo}.

\subsection{$\Gamma$-convergence and convergence of minimizers}

Statements (ii) and (iii) in Theorem \ref{beta2} can be summarized in terms of {\em $\Gamma$-convergence} \cite{dalmaso}, which is one of the basic notions of convergence in Calculus of Variations.
A sequence of functionals  $\{F_n:Z\to \bar{\mathbb{R}}\}_{n\to \infty}$ defined on a metric space $Z$, is said to $\Gamma$-converge to ${F}:Z\to \bar{\mathbb{R}}$ when
the following two conditions hold:
\begin{enumerate}[leftmargin=7mm]
\item[(i)] For any converging sequence $\{z_n\}_{n\to \infty}$ in $Z$ we have:
${F}\big(\lim_{n\to \infty} z_n\big) \leq \liminf_{n\to 0} {F}_n(z_n).$
\item[(ii)] For every $z\in Z$ there exists $\{z_n\}_{n\to \infty}$ converging
to $z$ and such that: ${F}(z) = \lim_{n\to \infty} {F}_n(z_n)$.
\end{enumerate}
We then write $F_n\Gc F$. It is an exercise to show that if, additionally, there exists a compact set $K\subset Z$ with the property that $\inf_Z F_n = \inf_K F_n$ for all $n$, then we have:
\begin{enumerate}[leftmargin=7mm]
\item[(i)] For any sequence $\{z_n\in K\}_{n\to \infty}$ of
  approximate minimizers to $F_n$, namely when $|F_n(z_n)-\inf_Z F_n|\to 0$, any accumulation point $z=\lim_{k\to\infty} z_{n_k}$ is a minimizer of $F$, i.e. $F(z) =\inf_X F$. In particular, $F$ has at least one minimizer and it has at least one minimizer in $K$.
\item[(ii)] For every minimizer $z\in Z$  of $F$, there
  exists a {\em recovery sequence} of approximate minimizers $z_n\to
  z$ so that $|F_n(z_n)-\inf_Z F_n|\to 0$.
\end{enumerate}


\medskip

In view of the compactness assertion (i), Theorem \ref{beta2} hence yields:
\begin{corollary}\label{convbeta2}
There holds, with respect to the strong convergence in $W^{1,2}(\Omega^1,\R^3)$:
$$\frac{1}{h^2}\mathcal{E}^h\big(y(x',hx_3)\big) \Gc 
\left\{\begin{array}{ll} \mathcal{I}_{2,g}(y)
& \mbox{ if $y\in W^{2,2}(\omega,\R^3)$ satisfies (\ref{iso2})} \\ +\infty & \mbox{ otherwise,} \end{array}\right. $$
Consequently, there is a one-to-one correspondence between limits of sequences of (global) approximate minimizers to the energies $\mathcal{E}^h$ and (global) minimizers of $\mathcal{I}_{2,g}$, provided that the induced metric $g(\cdot,0)_{2\times 2}$ has a $W^{2,2}$ isometric immersion from $\omega$ to $\R^3$.
\end{corollary}

It is useful to make a couple of observations. First, we point out that, in general, one cannot expect $\mathcal{E}^h$ to posses a minimizer. The lowersemicontinuity of the energy $\mathcal{E}$ in (\ref{nonE}) allowing for the direct method of Calculus of Variations, is tied to the quasiconvexity of the energy
density, whereas $F\mapsto\mbox{dist}^2(F, SO(3))$ is not even rank-one convex \cite[proof of Proposition 1.6]{zhang}.

\begin{figure}[htbp]
\centering
    \includegraphics[width=\textwidth]{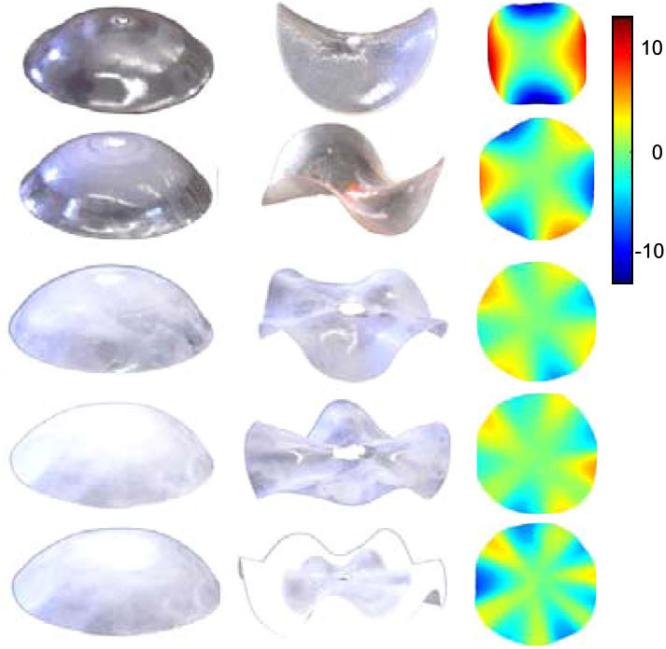}
    \vskip -.5in
\caption{As a function of the thickness of a swollen sheet, one can achieve either elliptical or hyperbolic geometries. In the limit of vanishing thickness, the shape converges to the limit implied by the $\Gamma$-convergence result in Theorem \ref{beta2} for the elliptic case (to a spherical dish), but shows an increasing preponderance to wrinkling on finer scales in the hyperbolic case. The multi-lobed swelling-induced wrinkling begs the question of the limiting behavior \cite[Reprinted with permission, copyright (2021) by the American Physical Society]{KVS}.}
\label{fig4}
\end{figure}

Second, we comment on the relation of Corollary \ref{convbeta2} with the experimental findings in \cite{KVS} that constructed a thickness-parametrised family of axially symmetric hydrogel disks (see Figure \ref{fig4}). The explicit control on the radial concentration $c(r)$ of the temperature-responsive polymer (N-isopropylacrylamide) resulted in the ability to control the (locally isotropic) shrinkage factors of distances $\eta(r) = \eta(c(r))$ and led to the target metric $g_{2\times 2}=\mathrm{d}r^2 + \kappa^{-1}\sin(\rho\kappa^{1/2})^2 \mathrm{d}\theta^2$ on the midplate $\omega=B(0,R)$, written in polar coordinates $(r,\theta)$ and in terms of the prescribed
constant Gaussian curvature $\kappa\equiv\pm 0.0011$. 
While decreasing the thickness $h$, all disks with $\kappa>0$ kept the same basic dome-like shape, with minor variations along the edge (see left column figures in Figure \ref{fig4}. The energy related to (\ref{Ih}) was observed to stabilize as $h\to 0$, approaching a constant multiple of $h^2$ and exhibiting equipartition between bending and stretching; Hence, discs with positive curvature minimize their energy via the scenario in Corollary \ref{convbeta2}, by settling near the isometric immersion that is of the lowest bending content.

On the other hand, for $\kappa<0$, the disks were observed to undergo a set of bifurcations in which the number of nodes (within a single wave configuration) increases and is roughly proportional to $h^{-1/2}$. Measuring the bending content in this case led to $\mathcal{E}^h\sim h$ which seems to be linked to a stretching-driven process: the sharp increase of the bending content is compensated by simultaneous decrease in the stretching content. Hence, hyperbolic disks minimize their energy via a set of bifurcations, despite the existence of the smooth immersions $y$. 

\begin{problem}
Analyze the possible origins of the diversity of behaviour of the elliptic and hyperbolic disks in \cite{KVS}, as well as the discrepancy between the experimentally observed linear in $h$ energy scaling and the quadratic scaling obtained in Corollary \ref{convbeta2}. The accuracy of the experimental determining the target metric $g$ is finite, and the sensitivity to perturbations seems to be more pronounced in the negative Gauss curvature regimes. 
\end{problem}

\subsection{Energy scaling from convex integration} \label{sec3.4}

A separate energy bound may be obtained by constructing deformations $u^h$ through the Kirchhoff-Love
extension of isometric immersions of regularity $\mathcal{C}^{1,\alpha}$. Existence of such is guaranteed 
by techniques of convex integration \cite{DIS} for all $\alpha<1/5$, and this threshold implies the particular energy scaling bound in Proposition \ref{pr_general_scaling} below.
If we could take $\alpha<1/3$ (corresponding to the so-called ``one step'' in each ``stage''
of the Nash-Kuiper iteration scheme), then the exponent would be $\beta<1$. If we could take $\alpha<1/2$ for the flexibility threshold as conjectured in \cite{DS_rev},
then $\beta<4/3$. Recall that existence of a $W^{2,2}$ isometric immersion implies that $\inf \Eh$ may be  further decreased to $Ch^2$.  

\begin{proposition}\label{pr_general_scaling} 
Assume that $\omega\subset\R^2$ is simply connected with
$\mathcal{C}^{1,1}$-regular boundary. Then:
$$\inf\Eh \leq Ch^\beta\qquad \mbox{for all } \; \beta<\frac{2}{3}.$$
 \end{proposition}
\begin{proof}
Fix $\alpha\in (0,1/5)$ and let $y\in\mathcal{C}^{1,\alpha}(\bar\omega,\R^3)$ 
satisfy (\ref{iso2}). Define the vector field
$\vec b_1 \in \mathcal{C}^{0,\alpha}(\bar \omega,\R^3)$ by (\ref{Cosserrat_vector}),
yielding the following auxiliary matrix field:
$$B = \left[\partial_1y, \; \partial_2y, \; \vec b \right] \in
\mathcal{C}^{0,\alpha}(\bar \omega,\R^{3\times 3})
\qquad\mbox{satisfying: } \; \det B>0, \; \;  B^TB=g(\cdot,0) \; \mbox{ in } \bar\omega.$$
The last assertion above implies that:
\begin{equation}\label{grr}
B(x')A(x', 0)^{-1}\in SO(3)\quad\mbox{ for all }\; x'\in\bar\omega.
\end{equation}
Regularize now $y, \vec b$ to $y_\epsilon, {b}_\epsilon\in\mathcal{C}^\infty(\bar\omega,\R^3)$
by means of the family of standard convolution kernels $\{\phi_\epsilon(x) =
\epsilon^{-2}\phi(x/\epsilon)\}_{\epsilon\to 0}$ where $\epsilon$ is a power of $h$ to be chosen later:
$$ y_\epsilon= y*\phi_\epsilon,\quad b_\epsilon=\vec b*\phi_\epsilon,
\quad B_\epsilon = B*\phi_\epsilon \quad \mbox{ and } \quad \epsilon=h^t.$$
We will utilize the following bound, resulting from the commutator estimate \cite[Lemma 1]{CDS}:
\begin{equation}\label{m1}
\begin{split}
\big\| B_\epsilon^TB_\epsilon - g(\cdot,0)\big\|_{\mathcal{C}^0(\omega)} & \leq \big\|
B_\epsilon^TB_\epsilon - (B^TB)*\phi_\epsilon \big\|_{\mathcal{C}^0(\omega)}  + 
\big\| g(\cdot,0)*\phi_\epsilon - g(\cdot,0)\big\|_{\mathcal{C}^0(\omega)}
\\ & \leq C\epsilon^{2\alpha} +C\epsilon^2\leq C\epsilon^{2\alpha},
\end{split}
\end{equation}
where the $C\epsilon^2$ bound results by Taylor expanding $g$ up to second order terms in $x_3$. Denoting $D_\epsilon=\big[\partial_1 b_\epsilon, \,\partial_2 b_\epsilon, \,0\big]$, we get the uniform bounds: 
\begin{equation}\label{m2}
\|B_\epsilon - B\|_{\mathcal{C}^0(\omega)} \leq C\epsilon^{\alpha},
\qquad \|D_\epsilon \|_{\mathcal{C}^0(\omega)} \leq C\epsilon^{\alpha-1}.
\end{equation}

Consider the sequence of deformations: $u^h\in \mathcal{C}^\infty(\bar \Omega^h,\R^3)$ in:
\begin{equation*}\label{recseq0}
u^h(x', x_3) = y_\epsilon(x') + x_3 b_\epsilon(x'), \quad \mbox{ so that: }\; \nabla u^h = B_\epsilon+ x_3 D_\epsilon.
\end{equation*}
In particular, $\|\nabla u^h(x', hx_3) - B(x')\|_{\mathcal{C}^0(\Omega^1)}
\leq C(\epsilon^{\alpha} +h\epsilon^{\alpha-1})$ and since:
$ A(x', hx_3)^{-1} = A(x')^{-1} + \mathcal{O}(h)$ for all $(x', x_3)\in\Omega^1$,
it follows by (\ref{grr}) that:
\begin{equation*} 
\begin{split}
\big\|\dist\big(\nabla &u^h A^{-1}, SO(3)\big)\big\|_{\mathcal{C}^0(\Omega^h)} \leq 
\big\|\nabla u^h A^{-1}(x', hx_3) - BA^{-1}(x',0)\big\|_{\mathcal{C}^0(\Omega^1)} 
\\ & \leq C(\epsilon ^\alpha + h\epsilon^{\alpha-1}+h) \to 0 \; \mbox{ as } h\to 0,
\end{split}
\end{equation*}
if only $h \epsilon^{\alpha-1}\to 0$. We then use polar
decomposition theorem and conclude that for some $R^h = R^h(x', hx_3)\in SO(3)$ there holds:
\begin{equation*} 
\begin{split}
R^h\nabla & u^h(x', hx_3) A(x', hx_3)^{-1} = \big(A^{-1}(\nabla
u^h)^{T}\nabla u^h A^{-1}\big)^{1/2}(x', hx_3) \\ & =
\big(A (x', hx_3) ^{-1}(B_\epsilon^T B_\epsilon(x') + \mathcal{O}(hD_\epsilon))A (x', hx_3) ^{-1}\big)^{1/2}
\\ & = \big((A(x',0)^{-1}+\mathcal{O}(h))(g(x',0) +
\mathcal{O}(\epsilon^{2\alpha} + h\epsilon^{\alpha-1})(A(x',0)^{-1}+\mathcal{O}(h))\big)^{1/2} 
\\ & = \big(Id_3 + \mathcal{O}(h+\epsilon^{2\alpha} + h\epsilon^{\alpha-1}\big)^{1/2}
=  Id_3 + \mathcal{O}(h+\epsilon^{2\alpha} + h\epsilon^{\alpha-1} ),
\end{split}
\end{equation*}
in virtue of (\ref{m1}) and (\ref{m2}). Consequently, we obtain the energy bound:
\begin{equation*}
\begin{split}
\Eh(u^h) & = \int_{\Omega^1}W\big(R^h \nabla
u^h(A^h)^{-1}(x', hx_3) \big)\;\mbox{d}(x', x_3) \\ & \leq
C \big(h + \epsilon^{2\alpha} + h\epsilon^{\alpha-1}\big)^2 = C \big(h + h^{2\alpha t} + h^{1+(\alpha-1)t}\big)^2.
\end{split}
\end{equation*}
Minimizing the right hand side above corresponds to
maximizing the minimal of the three displayed exponents. We hence choose $t$ in $\epsilon=h^t$ so that
${2\alpha t}= {1+(\alpha-1)t}$, namely $t=\frac{1}{\alpha+1}$. This leads to the estimate:
$$\inf \Eh \leq Ch^{\frac{4\alpha}{\alpha+1}}\quad\mbox{ for all }\; \alpha<\frac{1}{5},$$
which completes the proof.
\end{proof}

\begin{problem}
Analyze the limiting behavior of minimizing deformations in the intermediate energy scaling regime $\inf\Eh\simeq Ch^\beta$ for $ \beta\in [2/3, 2)$. Is it necessarily guided by an isometric immersion of some prescribed regularity? Find the $\Gamma$-limits of scaled energies $\frac{1}{h^\beta}\Eh$. 
\end{problem}

\subsection{The intermediate scalings}

As a point of comparison, we remark that higher energy scalings $\inf\mathcal{E}^h\sim h^\beta$ may result due to the sheet being forced at a boundary, due to the presence of external forces associated with gravity, the presence of an elastic substrate etc. all of which can lead to a range of microstructural patterns that are wrinkle-like. From Theorem \ref{beta2}, we recall that the regime $\beta \geq 2$ pertains to the ``no wrinkling'' family of almost minimizing deformations, that are perturbations of a $W^{2,2}$ isometric immersion. While the systematic description of the singular limits associated with exponents  $\beta<2$ is not yet available, there are a number of examples of the variety of emerging patterns that are illustrative. 

When a thin film is either clamped or weakly adhered to a substrate and subject to thermal or mechanical loads it can buckle and blister  \cite{JS, KoBe, BCDM, BCDM2}; in these cases, the energy scaling estimates yield $\beta=1$. A similar exponent is also seen in cases when a thin film wrinkles in response to metric constraints \cite{KoBe}, or forms a hanging drape exhibiting fluted patterns that coarsen as a function of distance from the point of support \cite{CMP, Bella}. In related experiments and theory,  when a thin shell of non-zero curvature is placed on a liquid bath, it forms complex wrinkling patterns \cite{TobascoARMA} with a range of $\beta$ between $0$ and $1$, depending on the strengths of the elastic and substrate forces. Moving from sheets or surfaces towards ribbons that have all three dimensions far from each other, papers
 \cite{Obrien, KO}  analyze wrinkling in the center of a stretched and twisted ribbon and find that $\beta=4/3$. Moving away from situations associated with non-local wrinkled microstructures, there have been a number of studies of localized structures associated with the theoretical and experimental analysis of {\em conical singularities} that arise in crumpled sheets \cite{CM2,CM3,CM4} that have been recently analyzed mathematically  \cite{Olber1, Olber2, Olber3} and lead to energetic estimates for this singularity of the form $\mathcal{E}^h\sim h^2\log(1/h)$.  And, in cases where the sheet is strongly creased, as in {\em origami patterns}, energy levels are associated with $\beta=5/3$ \cite{Maggi, Venka}. 
We remark that the mentioned papers do not address the dimension reduction, but rather analyze the chosen actual configuration of the prestrained sheet.

Closely related is also the literature on shape selection in non-Euclidean plates, exhibiting hierarchical {\em buckling patterns} in the limit of zero strain plates  with $\beta=2$, where the complex morphology is due to non-smooth energy minimization \cite{Ge1, Ge2, GSSV}.
Various geometrically nonlinear thin plate theories have been used to analyze the self-similar structures with metric asymptotically flat at infinity \cite{1b} that include a disk with edge-localized growth \cite{ESK1}, the shape of a long leaf \cite{18b}, or torn plastic sheets \cite{22a}.


\section{Hierarchy of limiting theories in the non-wrinkling regimes}\label{sec4}

\begin{figure}[htbp]
\centering
    \includegraphics[width=\textwidth]{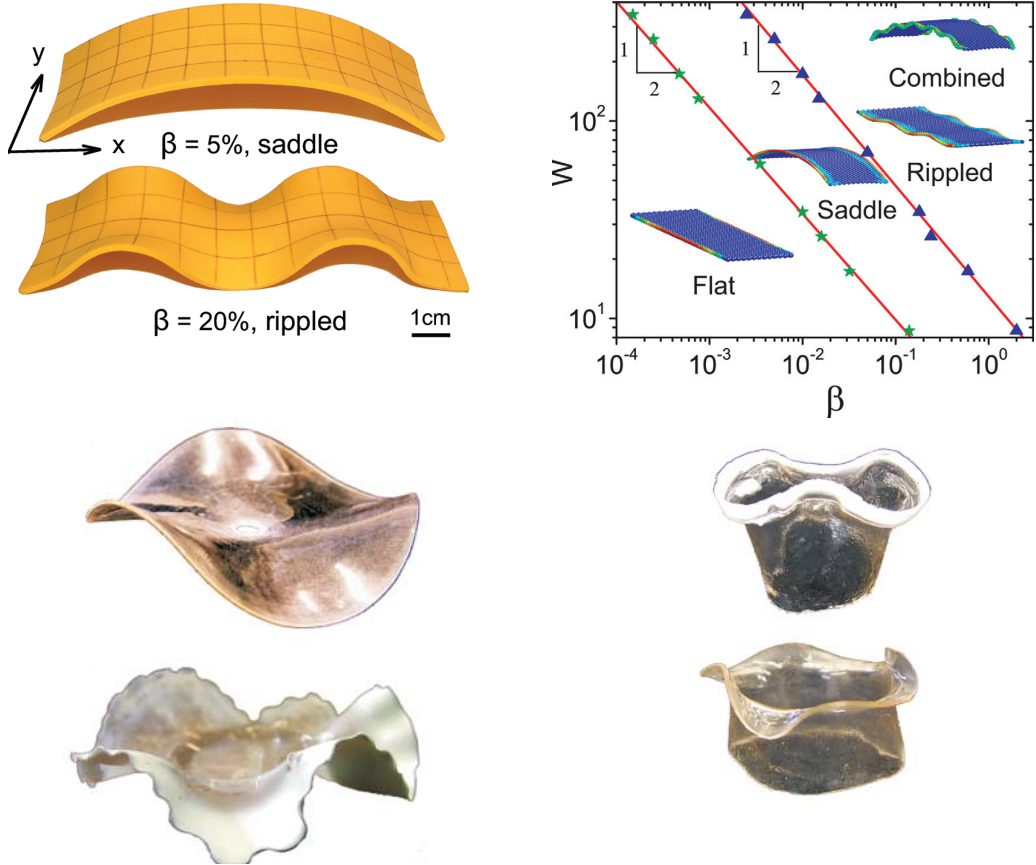}
\caption{A range of patterns arise when a thin sheet is inhomogeneously stretched plastically or swells in response to a solvent. (a) By dragging one's nails along the edges of a foam strip weakly, a flat surface transitions to one that is hyperbolic. (b) The same process carried out strongly leads to a surface that is strongly rippled, much like the edges of a leaf \cite[Copyright (2021) National Academy of Sciences]{18b}. (c) Thin sheets of a circular gel disk deform into a hyperbolic surface with two lobes. (d) Thinner sheets deform into multi-lobed sheets which are able to relieve the swelling-induced frustration by changing their curvature on multiple scales \cite[Reprinted with permission from AAAS]{27}.}
\label{fig5}
\end{figure}

We now detail the complete set of results relating the context of dimension reduction in
non-Euclidean elasticity with the quantitative immersability of Riemann metrics. As shown in Figure \ref{fig5}, a range of distinct behaviours of a thin sheet takes place in response to the prestrains of different orders. Within the formalism of finite elasticity, such patterns result from the sheet buckling to relieve growth or swelling induced by the residual strains. These will be measured via the scaling of the prestrain metric's Riemann curvatures, 
as explained below.

\subsection{The energy scaling quantization}

Observe that in view of Theorem \ref{beta2}, there holds $\lim_{h\to 0}\frac{1}{h^2}\inf\Eh = 0$
if and only if there exists $y\in W^{2,2}(\omega, \R^3)$ such that with $\vec b$ as in (\ref{Cosserrat_vector}):
\begin{equation}\label{compa0}
(\nabla y)^T\nabla y = g(\cdot, 0)_{2\times 2} \quad \mbox{ and }\quad
\mbox{sym}\big((\nabla y)^T\nabla \vec b\big) =
\frac{1}{2}\partial_3g(\cdot, 0)_{2\times 2} \quad \mbox{ in }\;\omega.
\end{equation}
The above compatibility of tensors $g(\cdot, 0)_{2\times 2}$
and $\partial_3g(\cdot, 0)_{2\times 2}$ is equivalent to the
satisfaction of the Gauss-Codazzi-Meinardi equations for the
related first and second fundamental forms:
$ I = (\nabla y)^T\nabla y$, $II= (\nabla y)^T\nabla\vec N =
\sqrt{g^{33}}\big(\mbox{sym}((\nabla y)^T\nabla \vec b) - \frac{1}{2}
\partial_3g(\cdot, 0)_{2\times 2} \big) - \frac{1}{\sqrt{g^{33}}}\big[\Gamma_{ij}^3(\cdot, 0)\big]_{i,j=1\ldots 2}$.
These three compatibility conditions turn out to be are precisely expressed by:
\begin{equation}\label{R12} 
\mathcal{R}_{12,12}(\cdot, 0) = \mathcal{R}_{12, 13}(\cdot, 0) = \mathcal{R}_{12, 23}(\cdot, 0) = 0\quad\mbox{ in }\;\omega.
\end{equation}

\begin{corollary}\cite{LRR,BLS,LL} \label{scal_1}
Condition $\frac{1}{h^2}\inf\Eh\to 0$ as $h\to 0$ is equivalent to: $\min\mathcal{I}_{2,g}=0$, and further to (\ref{R12}). In case (\ref{R12}) holds, we have: $\mbox{Ker}\;\mathcal{I}_{2,g}=\big\{Ry_0+c;~ R\in SO(3),
~c\in\R^3\big\}$, where $y_0:\bar \omega\to\R^3$ is a unique ``compatible'' smooth isometric immersion of $g(\cdot, 0)_{2\times 2}$ satisfying (\ref{compa0}) together with its corresponding Cosserat vector $\vec b = \vec b_1$. 
Moreover: $\inf\Eh\leq Ch^4$.
\end{corollary}

To justify the last assertion, we define the family $\{u^h\}_{h\to 0}$ as in the proof of Proposition \ref{pr_general_scaling}:
$$u^h(x', x_3) = y_0(x')+x_3\vec b_1(x') + \frac{x_3^2}{2}\vec b(x'), \quad\mbox{ so that: }\; \nabla u^h = B_0 + x_3B_1+\mathcal{O}(x_3^2).$$
By polar decomposition, the tensor $(\nabla u^h)A^{-1}$ coincides with $Z\doteq\big(\big((\nabla u^h)A^{-1}\big)^T \big((\nabla u^h)A^{-1}\big)\big)^{1/2}$ up to a rotation. Since:
\begin{equation*}
\begin{split}
& (\nabla u^h)^T\nabla u^h = B_0^TB_0 +2x_3\mbox{sym}(B_0^TB_1) + \mathcal{O}(x_3^2),\\
& A(x', x_3)^{-1}=A(x',0)^{-1}- x_3A^{-1}(\partial_3A)A^{-1}(x',0) + \mathcal{O}(x_3^2),
\end{split}
\end{equation*}
it follows that $Z^2$ equals:
\begin{equation*}
\begin{split}
 A^{-1}(\nabla u^h)^T&(\nabla u^h)A^{-1} = A^{-1}\big(g(x',0) +
2x_3\mbox{sym}(B_0^TB_1)+\mathcal{O}(x_3^2)\big)A^{-1} \\ & = Id_3 + 2x_3
A^{-1}\mbox{sym}\big( B_0^TB_1 - A\partial_3A\big) A^{-1}(x',0) +
\mathcal{O}(x_3^2) = Id_3 + \mathcal{O}(x_3^2).
\end{split}
\end{equation*}
The last equality above is achieved by choosing
$\vec b_2:\bar \omega\to \R^3$ such that $\mbox{sym}\big(B_0^TB_1 -
A\partial_3A (\cdot,0)\big)=0$, because the $2\times 2$ minor of the
indicated tensor is zero due to (\ref{compa0}). Consequently:
$$\Eh(u^h) = \frac{1}{h}\int_{\Omega^h}W(Z) \leq \frac{1}{h}\int_{\Omega^h}W\big(Id_3+
\mathcal{O}(x_3^2)\big) \leq Ch^4.$$

\smallskip

The following general result proves that 
the only viable scalings of $\inf \Eh\sim h^\beta$ in the regime
$\beta\geq 2$ are the even powers $\beta=2n$.

\begin{theorem}\label{th_quanti}\cite{Lew_last}
For every $n\geq 2$, if  $\lim_{h\to 0}\frac{1}{h^{2n}}\inf \Eh=0$ then $\inf \Eh\leq Ch^{2(n+1)}$.
Further, the following three statements are equivalent:
\begin{enumerate}[leftmargin=7mm]
\item[(i)] $\inf \Eh\leq Ch^{2(n+1)}$.
\item[(ii)] $\mathcal{R}_{12,12}(\cdot,0) = \mathcal{R}_{12,13}(\cdot,0) = \mathcal{R}_{12,23}(\cdot,0)=0$ and $\partial_3^{(k)} \mathcal{R}_{i3,j3}(\cdot,0) = 0$ in $\omega$,  
  for all $k=0\ldots n-2$ and all $i,j=1\ldots 2$. 
\item[(iii)] There exist smooth fields $y_0, \{\vec b_k\}_{k=1}^{n+1}:\bar\omega\to\mathbb{R}^3$ and
frames $\big\{B_k  = \big[ \partial_1 \vec b_k, ~\partial_2 \vec b_k, ~\vec
b_{k+1}\big]\big\}_{k=1}^n$, $B_0=\big[ \partial_1 y_0, ~\partial_2 y_0, ~\vec b_{1}\big]$, 
such that: $\displaystyle \sum_{k=0}^{m}\binom{m}{k} B_k^T B_{m-k}
-\partial_3^{(m)}g(\cdot,0) = 0\;$ for all $ m=0\ldots n$. Equivalently:
$\displaystyle \Big(\sum_{k=0}^n\frac{x_3^k}{k!}B_k\Big)^T\Big(\sum_{k=0}^n\frac{x_3^k}{k!}B_k\Big)=
g(x',x_3) + \mathcal{O}(h^{n+1})$ on $\Omega^h$ as $h\to 0$.
The field $y_0$ is the unique smooth isometric immersion of $g(\cdot,0)_{2\times 2}$ into
$\R^3$ for which $\mathcal{I}_{2,g}(y_0)=0$.
\end{enumerate}
\end{theorem}

We note that if $\mathcal{R}(\cdot, 0)= 0$ and $\partial_3^{(m)}
\big[\mathcal{R}_{i3,j3}(\cdot ,0)\big]_{i,j=1\ldots 2} = 0$  on $\omega$ for  all $m=0\ldots n-2$, but 
$\partial_3^{(n-1)} \big[\mathcal{R}_{i3,j3}(\cdot ,0)\big]_{i,j=1\ldots 2} \neq 0$, then: 
$ ch^{2(n+1)}\leq \inf \mathcal{E}^h\leq Ch^{2(n+1)}$ for  some $c,C>0.$
The conformal metrics $g(x', x_3) = e^{2\phi(x_3)}Id_3$ provide a class of examples for the viability of
all scalings: $\inf \Eh\sim h^{2n}$ by choosing $\phi^{(k)}(0)=0$ for
$k=1\ldots n-1$ and $\phi^{(n)}(0)\neq 0$. 

\bigskip

\begin{figure}[htbp]
\centering
\renewcommand{\arraystretch}{1.5}
\hspace{-6mm} $\begin{array}{l | l | l | l}
 {\beta} & ~~~~~\mbox{\begin{minipage}{3cm}  {asymptotic \\ expansion} \end{minipage}} &
~~~~~~\mbox{\begin{minipage}{4cm} {constraint / regularity}\end{minipage}} &
\mbox{limiting energy } {\mathcal{I}}_{\beta,g} \\ [10pt]\hline
 {2} & \hspace{-1.5mm}\begin{array}{ll} y(x') \\ \scriptstyle{\hspace{-1mm} 
\big\{ 3d:~ y(x') + x_3\vec b(x') \big\} } \end{array} & \hspace{-2mm} \footnotesize{\begin{array}{l} y\in W^{2,2}
\\ (\nabla y)^T\nabla y = g(x',0)_{2\times
  2} \end{array} } & \hspace{-2mm} \footnotesize{\begin{array}{l} c\|(\nabla
y)^T\nabla \vec b - \frac{1}{2}\partial_3g(x',0)_{2\times
  2}\|_{\mathcal{Q}_2}^2 \vspace{0.5mm} \\ \big[\partial_3y,\partial_2y,\vec b\big]\in
SO(3)g(x',0)^{1/2}\end{array} }\\ [15pt]\hline {4} & \hspace{-2mm}\begin{array}{l} y_0(x')+hV(x') \\
  \hspace{10.5mm} +\; h^2w^h(x')\end{array} &  \hspace{-2mm} \footnotesize{\begin{array}{l}
  \mathcal{R}_{12,12}, \mathcal{R}_{12,13}, \mathcal{R}_{12,23}(x',0) =0 \hspace{-2mm} \\ \big((\nabla y_0)^T\nabla
  V\big)_{sym}=0, \\  \big((\nabla y_0)^T\nabla w^h\big)_{sym}\to 
  \mathbb{S} \\ V\in W^{2,2}(\omega, \R), ~w^h\in W^{1,2}(\omega, \R^3) \end{array} } \hspace{-2mm} 
& \hspace{-2mm} \footnotesize{\begin{array}{l}
  c_1\|\frac{1}{2}(\nabla V)^T\nabla V + \mathbb{S} + \frac{1}{24}(\nabla \vec b_1)^T\nabla\vec b_1
  \\ \hspace{4mm} -\frac{1}{48}\partial_{33}g(x',0)_{2\times 2}\|_{\mathcal{Q}_2}^2 \\ +
  c_2\|(\nabla y_0)^T\nabla \vec p + (\nabla V)^T\nabla \vec b_1\|_{\mathcal{Q}_2}^2 \\ +c_3
  \|\big[\mathcal{R}_{i3,j3}(x',0)\big]_{i,j=1,2}\|_{\mathcal{Q}_2}^2 \vspace{2mm} \end{array}}
\\ [15pt] \hline  
\begin{array}{c} \hspace{-2mm} {6} \\
 \hspace{-2mm}  {\vdots} \end{array} & y_0(x') + h^2V(x')
&\hspace{-2mm} \footnotesize{\begin{array}{l} \mathcal{R}_{ab,cd}(x',0)=0\\ 
\big((\nabla y_0)^T\nabla V\big)_{sym}=0, ~ V\in W^{2,2} \hspace{-2mm}
\end{array} } & \hspace{-2mm} \footnotesize{\begin{array}{l}
  c_2\|(\nabla y_0)^T\nabla \vec p + (\nabla V)^T\nabla \vec
  b_1 +\alpha \big[\partial_3 \mathcal{R}\big] \|_{\mathcal{Q}_2}^2
\\ + c_3 \|\mathbb{P}_{\mathcal{S}_{y_0}^\perp}\big[\partial_3 \mathcal{R}\big]\|_{\mathcal{Q}_2}^2
+ c_4 \|\mathbb{P}_{\mathcal{S}_{y_0}}\big[\partial_3 \mathcal{R}\big]\|_{\mathcal{Q}_2}^2 
\end{array} } \\ [15pt]\hline \hspace{-1.5mm} 
\begin{array}{c} \hspace{-2mm} {2n} \\
 \hspace{-2mm}  {\vdots} \end{array} \hspace{-1.5mm} 
& \begin{array}{l} \hspace{-2mm} y_0(x')+h^{n-1}V(x') \\ \hspace{-3mm}
\scriptstyle{\big\{ 3d:~ y_0 + \sum_{k=1}^{n-1} \frac{x_3^k}{k!}\vec
  b_k(x') }\\ \hspace{4mm} \scriptstyle{+ h^{n-1}V(x') } \\ \hspace{4mm} \scriptstyle{+
  h^{n-1} x_3\vec p(x')\big\} } \end{array} \hspace{-4mm} 
& \hspace{-2mm} \footnotesize{ \begin{array}{l} \mathcal{R}_{ab,cd}(x',0)=0 \\ 
  \big[\partial_3^{(k)}\mathcal{R}\big](x',0) = 0 \;\forall k\leq n-3 \hspace{-2mm} \\ 
\big((\nabla y_0)^T\nabla V\big)_{sym}=0, ~ V\in W^{2,2}\end{array} \hspace{-2mm}}
& \hspace{-2mm} \footnotesize{ \begin{array}{l}
  c_2\|(\nabla y_0)^T\nabla \vec p + (\nabla V)^T\nabla \vec b_1
+ \alpha \big[\partial_3^{(n-2)} \mathcal{R}\big] \|_{\mathcal{Q}_2}^2
\\ + c_3 \|\mathbb{P}_{\mathcal{S}_{y_0}^\perp}\big[\partial_3^{(n-2)}\mathcal{R}\big]\|_{\mathcal{Q}_2}^2 \\
+ c_4 \|\mathbb{P}_{\mathcal{S}_{y_0}}\big[\partial_3^{(n-2)}R\big]\|_{\mathcal{Q}_2}^2 
\end{array} } \end{array}$ 
\caption{The infinite hierarchy of $\Gamma$-limits for prestrained
  films, scaling $\beta\geq 2$}
\label{table_prestrain}
\end{figure}

\medskip

\subsection{The infinite hierarchy of $\Gamma$-limits}\label{sec_compa2} 

To derive the counterpart of Corollary \ref{convbeta2} for higher energy scalings, one observes the following compactness properties under the assumption that $\Eh(u^h) \le Ch^{2(n+1)}$, for some $n\geq 1$.
First, \cite{Lew_last}, there exist $c^h\in \mathbb R^3$, $R^h\in SO(3)$ such that:
$$\displaystyle V^h(x') = \frac{1}{h^n}\fint_{-h/2}^{h/2} (\bar
R^h)^T\big(u^h(x', x_3) - c^h\big) - \Big(y_0(x') +
\sum_{k=1}^n\frac{x_3^k}{k!}\vec b_k(x')\Big)~\mbox{d}x_3$$
converge as $h\to 0$  in  $W^{1,2}(\omega, \R^3)$, to a limiting
displacement $V$ that is an infinitesimal isometry:
$$V\in\mathcal{V}_{y_0}= \big\{V\in W^{2,2}(\omega, \R^3);~ \mbox{sym}\big((\nabla
y_0)^T\nabla V\big) =0\big\}.$$ 
In particular, there exists $\vec p\in W^{1,2}(\omega, \R^3)$ with 
$\mbox{sym}\big(B_0^T\big[\nabla V, ~ \vec p\big]\big) = 0$. 
Second, the strains:
$$\frac{1}{h}\mbox{sym}\big((\nabla y_0)^T\nabla V^h\big)$$
converge as $h\to 0$, weakly in $L^2(\omega, \R^{2\times 2})$ to a limiting $\mathbb{S}$ in the finite strain space:
$$ \mathbb{S}\in\mathcal{S}_{y_0}= \mbox{cl}_{L^2}\big\{
\mbox{sym}((\nabla y_0)^T\nabla w_n);~ w_n\in W^{1,2}(\omega, \R^3)\big\}.$$ 
The space $\mathcal{S}_{y_0}$ can be
  identified, in particular,  in the following two cases on $\omega$
  simply connected. When $y_0=id_2$, then $\mathcal{S}_{y_0}= \{\mathbb{S}\in
L^2(\omega, \mathbb{R}^{2\times 2}_{sym}); ~\mbox{curl}\,\mbox{curl}\,\mathbb{S}=0\}$.
When the Gauss curvature $\kappa((\nabla y_0)^T\nabla
y_0)=\kappa\big(g(\cdot, 0)_{2\times 2})>0$ on $\bar
\omega$, then $\mathcal{S}_{y_0}= L^2(\omega, \mathbb{R}^{2\times
  2}_{sym})$, as shown in \cite{lemopa_convex}. 

\smallskip

We further have the following $\Gamma$-convergence results with respect to the above compactness statements. The infinite hierarchy of the limiting prestrained theories is gathered in Figure \ref{table_prestrain}.

\begin{theorem}\cite{LL, Lew_last}
In the energy (\ref{Ih}) scaling regimes indicated in
Theorem \ref{th_quanti}, the following holds. For the von K\'arm\'an-like
regime, we have for all $V\in\mathcal{V}_{y_0}$ and $\mathbb{S}\in\mathcal{S}_{y_0}$:
\begin{equation*}
\begin{split}
\frac{1}{h^4}\Eh & \Gc \mathcal{I}_{4,g}(V,\mathbb{S}) \\ = &  \, \frac{1}{2}\int_{\omega}\mathcal{Q}_2\Big(x', 
\underbrace{\mathbb{S}(x') + \frac{1}{2}\nabla V(x')^T\nabla V(x') +
  \frac{1}{24} \nabla \vec b_1(x')^T\nabla\vec b_1(x')-
  \frac{1}{48} \partial_{33}g(x', 0)_{2\times 2}}_{\mbox{stretching}}\Big) ~\mathrm{d}x' \\ & 
+ \frac{1}{24}\int_{\omega}\mathcal{Q}_2\Big(x', \underbrace{\nabla y_0(x')^T\nabla \vec p(x')+ \nabla
V(x')^T\nabla\vec b_1(x')}_{\mbox{bending}} \Big)~\mathrm{d}x' \\ & 
+\,\frac{1}{1440}\int_{\omega} \mathcal{Q}_{2} \Big(x',
  \underbrace{\left[\begin{array}{cc}\mathcal{R}_{13,13} & \mathcal{R}_{13,23}\\ \mathcal{R}_{13,23}
        & \mathcal{R}_{23,23}\end{array}\right]}_{\mbox{curvature}} \Big) ~\mathrm{d}x'. 
\end{split}
\end{equation*}
For all $n\geq 1$ (which is the case parallel to linear elasticity),
we have for all $V\in\mathcal{V}_{y_0}$:
\begin{equation*}
\begin{split}
\frac{1}{h^{2(n+1)}}\Eh & \Gc \mathcal{I}_{2(n+1),g}(V)  \\ = & \, \frac{1}{24}\int_\omega
\mathcal{Q}_2\Big(x', \underbrace{(\nabla y_0)^T\nabla \vec p + (\nabla V)^T\nabla \vec b_1 + \alpha_n
\big[\partial_3^{(n-1)} \mathcal{R}_{i3,j3}\big]_{i,j=1\ldots
  2}}_{\mbox{bending}} \Big)~\mathrm{d}x'\qquad \\ & + \beta_n 
\int_\omega\mathcal{Q}_2\Big(x', \mathbb{P}_{\mathcal{S}_{y_0}^{\perp}}
\big(\big[\partial_3^{(n-1)} \mathcal{R}_{i3,j3}\big]_{i,j=1\ldots  2}\big)\Big) ~\mathrm{d}x' \\ & + \gamma_n 
\int_\omega\mathcal{Q}_2\Big(x', \mathbb{P}_{\mathcal{S}_{y_0}}
\big(\big[\partial_3^{(n-1)} \mathcal{R}_{i3,j3}\big]_{i,j=1\ldots 2}\big) \Big) ~\mathrm{d}x',
\end{split}
\end{equation*}
where $\mathbb{P}_{\mathcal{S}_{y_0}}$, $\mathbb{P}_{\mathcal{S}_{y_0}^{\perp}}$ denote
orthogonal projections onto $\mathcal{S}_{y_0}$ and onto its $L^2$-orthogonal complement $\mathcal{S}_{y_0}^{\perp}$. The coefficients $\alpha_n, \beta_n, \gamma_n\geq 0$ are given explicitly and $\alpha_n\neq 0$ if and only if $n$ is even.
\end{theorem}

The functional $\mathcal{I}_{4,g}$ is indeed a von K\'arm\'an-like energy, consisting of stretching and bending (with respect to the unique, up to rigid motions, smooth isometric immersion $y_0$ that has zero energy in the prior $\Gamma$-limit (\ref{I2})) plus a new term quantifying the remaining three Riemann curvatures.
When $g=Id_3$ then $\mathcal{I}_{4,g}(V,\mathbb{S})$ reduces to the
classical von K\'arm\'an functional, given in terms of the
out-of-plane displacement $v$ in $V=(\alpha x^\perp+\beta,v)$ for which $\vec p = (-\nabla v, 0)$, and
the in-plane displacement $w$ in $\mathbb{S}=\mbox{sym}\nabla w$:
\begin{equation}\label{vonKar} 
\mathcal{I}_4(v,w)=\frac{1}{2}\int_\omega\mathcal{Q}_2\big(
\mbox{sym} \nabla w + \frac{1}{2}\nabla v\otimes \nabla
v\big)~\mbox{d}x' + \frac{1}{24}\int_\omega\mathcal{Q}_2(\nabla^2v)~\mbox{d}x'.
\end{equation}
We point out in passing that in \cite{DB, DCB}, a variant of the F\"oppl-von K\'arm\'an equilibrium equations has been formally derived from finite incompressible elasticity, via the multiplicative decomposition of deformation gradient \cite{RHM} used in finite plasticity \cite{Lee} and hyperelastic growth.

Likewise, each $\mathcal{I}_{2(n+1), Id_3}$ reduces to the classical linear elasticity: 
\begin{equation}\label{linlin}
\mathcal{I}_{2(n+1)}(v)= \frac{1}{24}\int_\omega \mathcal{Q}_2\big(\nabla^2 v\big)~\mbox{d} x'.
\end{equation}

In the present geometric context, the bending term $(\nabla
y_0)^T\nabla \vec p + (\nabla V)^T\nabla \vec b_1$ in $\mathcal{I}_{2(n+1),g}$ is of order $h^{n}x_3$ 
and it interacts with the curvature $\big[\partial_3^{(n-1)}
\mathcal{R}_{i3,j3}(\cdot,0)\big]_{i,j=1\ldots 2}$ which is of order $x_3^{n+1}$.
The interaction occurs only when the two terms have the same parity in
$x_3$, namely at even $n$, so that $\alpha_n=0$ for all $n$ odd. 
The two remaining terms measure the 
$L^2$ norm of $\big[\partial_3^{(n-1)} R_{i3,j3}(\cdot,0)\big]_{i,j=1\ldots 2}$, with distinct weights
assigned to  $\mathcal{S}_{y_0}$ and
$\big(\mathcal{S}_{y_0}\big)^\perp$ projections, again according to the parity of $n$. 
We also have:
$\inf_{\mathcal{V}_{y_0}}\mathcal{I}_{2(n+1),g} \sim \big\|\big[\partial_3^{(n-1)} \mathcal{R}_{i3,j3}(\cdot,0)\big]_{i,j=1\ldots 2}\|_{L^2(\omega)}^2.$

\begin{remark}
Parallel general results can be derived in the abstract setting of Riemannian manifolds: in \cite{KS14, KM14} $\Gamma$-convergence statements were proved for any dimension ambient manifold and codimension midplate, in the scaling regimes 
$\mathcal{O}(h^2)$ and $\mathcal{O}(1)$, respectively. In \cite{MS18}, the authors analyze scaling orders $o(h^2)$, $\mathcal{O}(h^4)$ and $o(h^4)$.
\end{remark}

\section{Floral morphogenesis, weak prestrain and special solutions of Monge-Amp\`ere equations}\label{sec-5}

\begin{figure}
    \centering
    \includegraphics[width=0.9\textwidth]{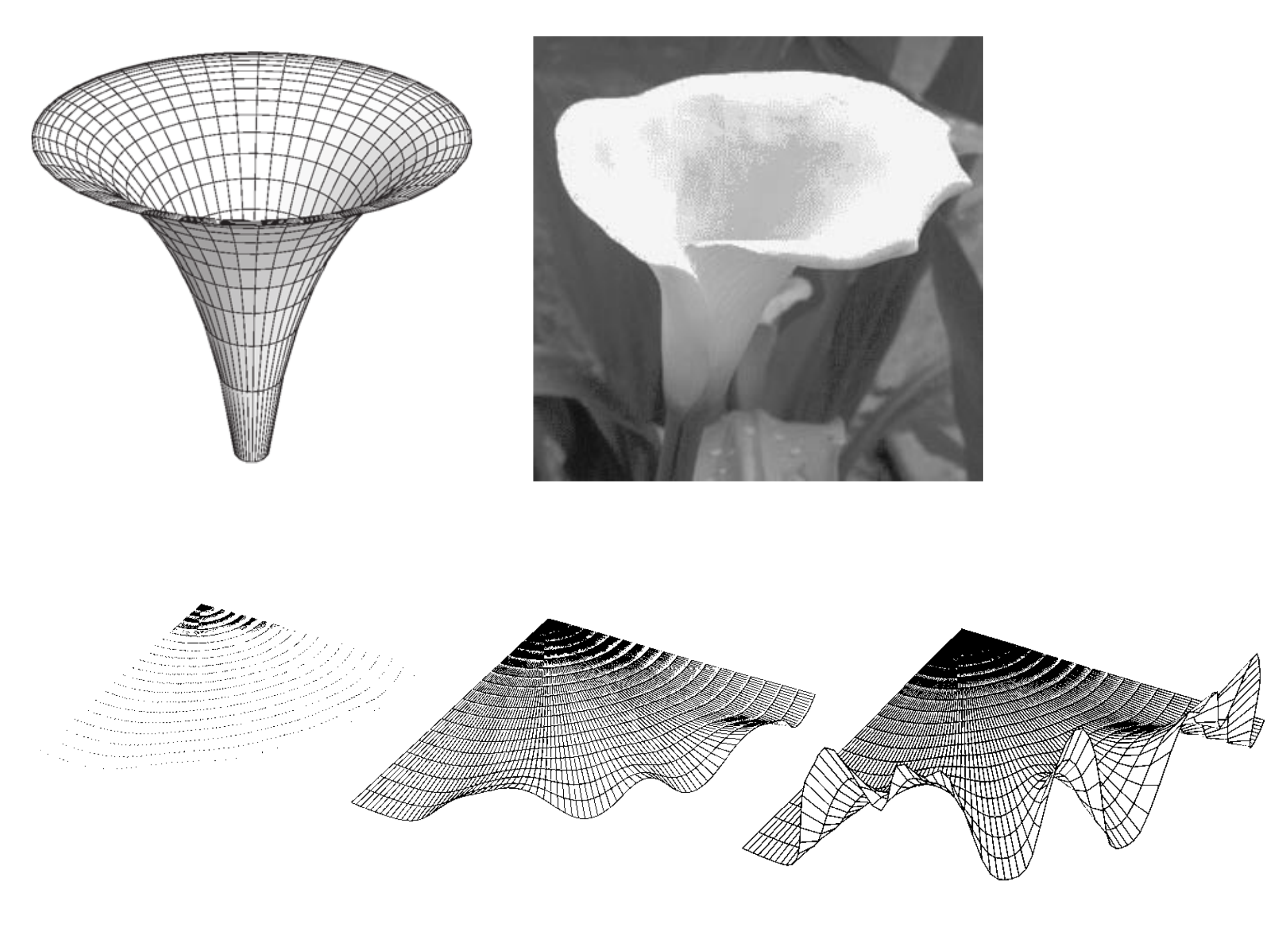}
    \vskip 0.in
    \caption{(a) Pseudosphere. (b) Picture of a calla lily. (c) Sample plots of the Jacobian function, encoding the Gauss curvature of the prestrain \cite{NV01}. }
    \label{nech_fig}
\end{figure}

We digress in this section to consider an interesting set of questions inspired by the remarkable examples of floral morphogenesis resembling parts of a pseudosphere (see Figure \ref{nech_fig}) altered by the presence of ripples along the free boundary. Early work \cite{NV01}, revisited in \cite{Nech_book} suggested that information on the profile of the boundary of a plant's leaf fluctuating in the direction transversal to the leaf's surface, can be read from the Jacobian of the conformal mapping corresponding to an isometric embedding of the given prestrain metric. This leads to the question of constructing solutions to the classical Monge-Amp\`ere equation, without prescribed boundary conditions but approximating the smallest bending content possible while preserving the regularity that allows for the consistent association of this bending content. 

A similar point of view has been adopted in \cite{GSSV}, for the choice of the target midplate metric: $g_{2\times 2}=Id_2 + 2\epsilon^2f(x_2)\mathrm{d}x_1^2$, posed on the infinite strip $\omega=\mathbb{R}\times [0,W]$. The coefficient field $f(x_2)$ corresponds to the $x_2$-dependent growth in the $x_1$ direction, localized near the $x_2=0$ edge of the sheet. An interesting class of buckling patterns that lower the bending energy of the sheet while satisfying the approximate isometry condition, was constructed via introducing ``branch point'' singularities, resulting in the multiple asymptotic directions, into solutions to the Gauss curvature constraint equation: $\det\nabla^2 v = \kappa(g_{2\times 2}) = -f''$. 

For weakly hyperbolic sheets with constant $\kappa<0$, the same construction has been recently refined in \cite{SV,Venka21}, using a discrete differential geometric approach linked with the notion of index of topological defects, to argue that the branch points are energetically preferred and may lead to the fractal-like recursive buckling patterns seen in some flowers and leaves.

\begin{problem}
While we will consider the problem here solely from a static elastic perspective, it is worth asking an allied question: how does a growing front leave behind a partially relaxed shape, i.e. that of a flower? 
\end{problem}

\subsection{Weak prestrain and the Monge-Amp\`ere constrained theories}

We assume that the given prestrain tensor $A^h=(g^h)^{1/2}$ on $\Omega^h$ is incompatible only through a perturbation of order which is a power of the film's thickness $h$:
\begin{equation}\label{ahform}
A^h(x', x_3)={Id}_3 +  h^\gamma S(x') + h^{\gamma/2}x_3 B(x').
\end{equation}
Here, $S,B:\bar\omega\to\R^{3\times 3}_{sym}$ are smooth tensors that correspond to stretching and bending with the choice\footnote{The more general choices of
  exponents $\alpha/2, \gamma/2$ were analyzed in \cite{JL}.}  of the exponents $\gamma, \frac{\gamma}{2}$. In this context, the counterpart of Theorem \ref{beta2} is as follows:

\begin{theorem}\cite{LOP}\label{th_Ah}
Let $u^h\in W^{1,2}(\Omega^h,\mathbb{R}^3$) satisfy: $\Eh (u^h) \leq Ch^{\gamma+2}$, for some $\gamma\in (1,2)$. 
\begin{enumerate}[leftmargin=7mm]
\item[(i)]  (Compactness). There exist $R^h\in SO(3)$, $c^h\in\mathbb{R}^3$ such that the following holds for 
$\{y^h(x',x_3) =  R^h u^h(x',hx_3) - c^h\}_{h\to 0}$. First,
$y^h$ converge to $x'$ in $W^{1,2}(\Omega^1,\mathbb{R}^3)$. 
Second, the scaled displacements:
${V^h(x')=\frac{1}{h^{\gamma/2}}\fint_{-1/2}^{1/2}y^h(x',t) -
  x'~\mathrm{d}t}\in W^{1,2}(\omega, \mathbb{R}^3)$
converge, up to a subsequence, to a displacement field $V$ of the form $V = (0,0,v)^T$, satisfying:
\begin{equation}\label{MA1}
v\in W^{2,2}(\omega,\mathbb{R}) \qquad \det\nabla^2 v = - \mathrm{curl}\, \mathrm{curl} \,S_{2\times 2}.
\end{equation} 

\item[(ii)] ($\Gamma$-convergence). If $\omega$ is simply connected with $\mathcal{C}^{1,1}$ boundary, then we have, with the same quadratic forms $\mathcal{Q}_2$ defined in (\ref{Qform}):
\begin{equation}\label{MA-bd}
\frac{1}{h^{\gamma+2}} \Eh(u^h) \Gc \mathcal{I}_{S, B}(v) = 
\frac{1}{12} \int_\omega \mathcal{Q}_2\big(x',\nabla^2 v + B_{2\times 2}\big)~\mathrm{d}x'. 
\end{equation}
\end{enumerate}
\end{theorem}

Similarly to Corollary \ref{isoim}, one can further deduce:

\begin{corollary}\label{isoMA}
The Monge-Amp\'ere problem (\ref{MA1}) has a solution $v\in W^{2,2}$ iff $\inf\Eh\leq C h^{\gamma+2}$.
Moreover, $c h^{\gamma+2} \leq \inf \Eh  \leq h^{\gamma+2}$ for some $c,C>0$ is equivalent to the solvability of (\ref{MA1}) and the simultaneous non-vanishing
of the lowest order terms (i.e. terms of order $\gamma$ and
$\frac{\gamma}{2}$, respectively) in $\mathcal{R}_{12,12}(\cdot, 0)$ and
$[\mathcal{R}_{12,i3}(\cdot, 0)]_{i=1,2}$. This last condition is equivalent to: 
$$\mathrm{curl}\,\mathrm{curl} \,S_{2\times 2} + \mathrm{det}\, B_{2\times 2} \not\equiv 0 \quad \mbox{ or } \quad \mathrm{curl}\, B_{2\times 2} \not\equiv 0 \quad \mbox{ in } \; \omega.$$
\end{corollary}

The equation in (\ref{MA1}) can be seen as an
equivalent condition for the family of deformations on $\omega$ (which, indeed,
corresponds to the recovery sequence in Theorem \ref{th_Ah} (ii)) given through the out-of-plane displacement
$v:\omega\to\R$, and any in-plane displacement $w:\omega\to\R^2$:
$$\phi^h(x') = \big(x' +h^\gamma w(x'), h^{\gamma/2}v(x')\big): \omega\to \R^3$$
to match the metric $g^h(\cdot,0)_{2\times 2}$ at the lowest order terms of its Gauss curvature. Indeed:
\begin{equation}\label{MA_exp}
\begin{split}
\kappa\big((\nabla \phi^h)^T\nabla \phi^h \big) & = \kappa\Big( Id_2+ h^\gamma (\nabla v\otimes \nabla v + 2\,\mbox{sym}\nabla w ) + h^{2\gamma}(\nabla w)^T\nabla w\Big) \\ & =
-\frac{h^{\gamma}}{2}\mbox{curl}\,\mbox{curl}\big(\nabla v\otimes \nabla v + 2\,\mbox{sym}\nabla w \big) + o(h^\gamma) = h^\gamma \det\nabla^2 v + o(h^\gamma), \\
\kappa\big(g^h(\cdot,0)_{2\times 2}\big) & = \kappa\Big(Id_2+
2h^\gamma S_{2\times 2} + h^{2\gamma}(S^2)_{2\times 2}\Big) = -h^{\gamma} \mbox{curl}\,\mbox{curl}\, S_{2\times 2} + o(h^{\gamma}).
\end{split}
\end{equation}
Recalling that the kernel of the operator {\em $\mbox{curl}\,\mbox{curl}$}
consists precisely of $\mbox{sym}\nabla w$, we further observe that (\ref{MA1}) is equivalent to the possibility of choosing $w$ such that $\phi^h$ is an isometric immersion of $(\omega, g^h(\cdot,0)_{2\times 2})$ at the leading order terms:
$$(\nabla\phi^h)^T\nabla \phi^h = Id_2+ 
2 h^\gamma \big(\frac{1}{2}\nabla v\otimes \nabla v + \mbox{sym}\nabla w\big) +\mathcal{O}(h^{2\gamma}) =  g^h(\cdot,0)_{2\times 2} + o(h^{\gamma}). $$

\smallskip

The above analysis suggests to view the Monge-Amp\'ere equation
$\det\nabla^2v = f$ through its very weak form, well defined
for all $v\in W^{1,2}_{loc}(\omega, \R)$, in the sense of distributions: 
\begin{equation} \label{MA}
{\mathcal{D}et}\, \nabla^2 v \doteq -\frac 12 \mbox{curl}\,\mbox{curl}  (\nabla v \otimes
\nabla v) =  f \qquad \mbox{in } \,\omega.
\end{equation}
Similarly to the results described in section \ref{sec3.4}, one can
then apply techniques of convex integration and show \cite{LP, DIS} that 
for any smooth $f :\bar \omega\to \R$ and $\alpha<\frac{1}{5}$, 
the set of $\mathcal{\mathcal{C}}^{1,\alpha}(\bar\omega)$ solutions to
(\ref{MA}) is dense in  $\mathcal{\mathcal{C}}^0(\bar\omega)$. 
That is,  for every $v_0\in \mathcal{\mathcal{C}}^0(\bar\omega)$ there exists a sequence
$v_n\in\mathcal{\mathcal{C}}^{1,\alpha}(\bar\omega)$, converging uniformly to $v_0$ and satisfying:
$\mathcal{D}et \nabla^2 v_n  = f$.
One consequence of this result is that the operator $\mathcal{D}et \nabla^2$
is weakly discontinuous everywhere in $W^{1,2}(\omega)$. By an
explicit construction, there follows a counterpart of Proposition \ref{pr_general_scaling}:

\begin{proposition}\label{pr_general_scaling2} \cite{JL}
Assume that $\omega\subset\R^2$ is simply connected with
$\mathcal{C}^{1,1}$ boundary. Then:
\begin{equation*}
\begin{split}
& \inf\Eh \leq Ch^\beta\qquad \mbox{for all } \; \gamma\in
\big[\frac{2}{7},2\big] \;\mbox{ and } \;\beta<\frac{5}{3}\gamma+\frac{2}{3},
\\ & \inf\Eh \leq Ch^\gamma \qquad \mbox{for all } \; \gamma\in \big(0, \frac{2}{7}\big).
\end{split}
\end{equation*}
\end{proposition}

\begin{problem}
Analyze the intermediate energy scaling regime $\inf\Eh\simeq h^\beta$ for $\beta\in \big[\frac{5}{3}\gamma+\frac{2}{3},\gamma+2\big)$ and
find the $\Gamma$-limits of the scaled energies $\frac{1}{h^\beta}\Eh$.
\end{problem}

\begin{problem}
Consider the generalization of  (\ref{MA}) to problems posed on higher-dimensional domains $\omega\subset\mathbb{R}^N$, in the context of the dimension reduction and isometry matching.
As shown in \cite{HL_MA}, the set $\{\mbox{sym}\nabla w;
~W^{1,2}(\omega,\mathbb{R}^N)\}$ is the kernel of the operator
$\it{Curl}^2$, where for $A\in L^2(\omega, \mathbb{R}^{N\times N})$ the $4$-tensor: $\it{Curl}^2(A)= \big[\it{Curl}^2(A)_{ab,
  cd}\big]_{a,b,c,d=1\ldots N}$ is given as the application of two exterior derivatives in: $\big[\partial_a\partial_cA_{bd} + \partial_b\partial_dA_{ac} - \partial_a\partial_dA_{bc} - 
\partial_b\partial_cA_{ad} \big]_{a,b,c,d}$.  Similarly to the calculation in (\ref{MA_exp}), there holds:
$\mathcal{R}_{ab,cd}(Id_N+\epsilon^2 A) = -\frac{\epsilon^2}{2}\it{Curl}^2(A)_{ab,
  cd} + o(\epsilon^2)$. Taking $A=\nabla v\otimes\nabla v$, one obtains that a scalar displacement field $v$ can be matched by a higher order perturbation vector field $w$, so that defining: $\bar\phi^h(x') =\big(
x'+h^2w (x'), hv(x')\big):\omega\to\mathbb{R}^{N+1}$, the given weak prestrain metric is matched by the pull-back metric in $(\nabla\bar\phi)^T\nabla\bar\phi = Id_N+h^2A+\mathcal{O}(h^4)$, if and only if:
$\big[\det(\nabla^2v)_{ab,cd}\big]_{a,b,c,d=1\ldots N}=-\it{Curl}^2(A)$.
\end{problem}

\subsection{Dimension reduction with transversely oscillatory prestrain}\label{dimosc}
We also mention the ``oscillatory setting'' where $g^h=(A^h)^2$ satisfy the following structure assumption: 
$$g^h(x', x_3) = \mathcal{G}^h(x', \frac{x_3}{h}) = \bar{\mathcal{G}}(x') + h \mathcal{G}_1(x',
\frac{x_3}{h})+ \frac{h^2}{2} \mathcal{G}_2(x', \frac{x_3}{h}) + \ldots\quad \mbox{ for all } ~x=(x',x_3)\in \Omega^h.$$ 
This set-up  includes the subcase $g^h=g$ of section \ref{sec4} upon taking:
$\bar{\mathcal{G}}_1= g(\cdot,0)$, $\mathcal{G}_1(x', t) = t\partial_3
g(x',0)$, $\mathcal{G}_2(x', t) = t^2\partial_{33}g(x',0)$, etc.
In \cite{LL} connections between these two cases were exhibited, via projections of appropriate
curvature forms on the polynomial tensor spaces and reduction to the
``effective non-oscillatory cases'' in the Kirchhoff-like ($h^2$) and von K\'arm\'an-like ($h^4$)
regimes. Compactness statements as in section \ref{sec_compa2}
are then still valid, with the
$\Gamma$-limits that consist of energies $\mathcal{I}_{2(n+1), \bar g}$ written for effective
metrics $\bar g$, plus the new ``excess term'' measuring the averaged deviation of $g^h$ from $\bar g$.

\begin{problem}
Derive the hierarchy of all the limiting theories in the oscillatory setting.
\end{problem}

\section{Classical geometrically nonlinear elasticity without prestrain}
\label{sec_clas}

\begin{figure}
    \centering
    \includegraphics[width=1.1\textwidth]{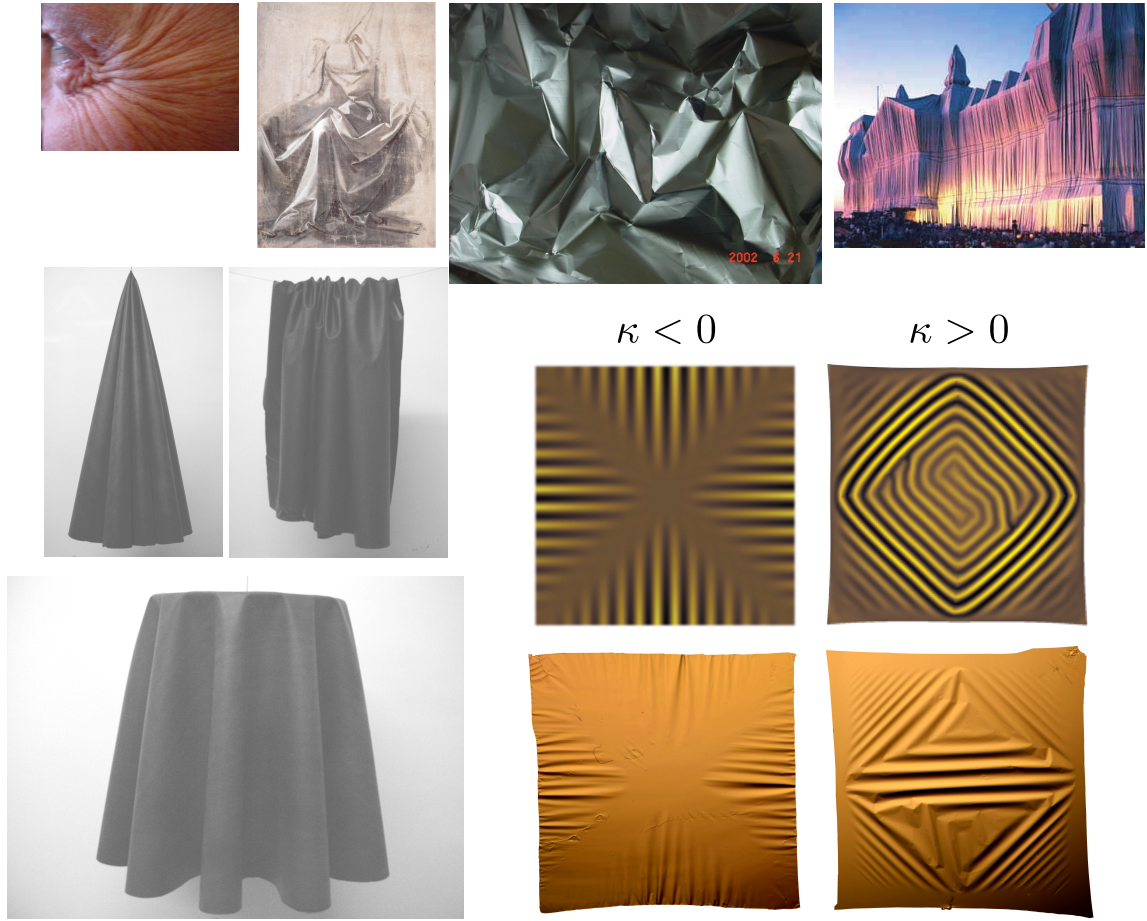}
    \vskip 0.in
    \caption{Wrinkles, drapes and crumples in thin sheets over a range of scales arise due to boundary and bulk forces. (a) Wrinkles in the neighborhood of the eye are driven by the muscular contractions. (b) The drape of a heavy piece of cloth on a knee is due to the combination of gravity and the presence of obstacles (a chiaroscuro by Leonardo). (c) The crumples in a sheet are reminiscent of the drape, but arise due to confinement, and are dominated by the presend  (d) The nearly uniform wrinkles on a fabric that wraps the Reichstag in Berlin (an inspiration of the artist Christo) are due to the presence of a series of horizontal ropes; otherwise the wrinkles will coalesce into larger and larger ones. (e) The elements of all drapes are a combination of the (in)ability to drape a point (such as tent-pole), a line (such as a curve) and a curve (such as a waist or a table), in the presence of gravity \cite[Copyright (2021) National Academy of Sciences]{CMP}. (f) Complex wrinkles also arise when non-Euclidean surfaces are flattened, as shown here for a patch of a surface that is either saddle-shaped ($\kappa < 0$) or spherical ($\kappa > 0$) in its natural configuration. Top: simulations; bottom: experiments. Images courtesy of Tobasco, Timounay, Todorova, Paulsen, and Katifori \cite{Tob}. }
    \label{fig6}
\end{figure}

When a thin plate or shell is constrained at the boundary, it can buckle, wrinkle or crumple depending on the nature and extent of the forcing. Similarly, when a plate or shell is subject to body forces such as those due to gravity in such contexts as draping a complex body, the sheet again folds and wrinkles in complex ways. Examples of the resulting patterns are shown in  Figure \ref{fig6}, and highlight the occurrence of three constituent building blocks: extended zones of short wavelength wrinkles, strongly localized conical structures, and the ridge-like structures that can arise either together or separately from the wrinkles. What is the hierarchy of limiting elastic theories in such situations? 

\subsection{The set-up and the finite hierarchy of $\Gamma$-limits for plates}
In this section we parallel the discussion of the hierarchy of the non-Euclidean thin films  presented in sections \ref{sec-3} - \ref{sec-5}.
Let $S\subset\mathbb{R}^3$ be a bounded, connected, oriented two-dimensional surface 
with unit normal $\vec n$. Consider a family $\{S^h\}_{h\to 0}$ of
thin shells around the midsurface $S$:
\begin{equation*}
S^h = \{x + t\vec n(x); ~~ x\in S, ~ -h/2< t < h/2\}.
\end{equation*}
The elastic energy (with density $W$ that satisfies (\ref{elastic_dens})) of a deformation $u^h: S^h\to\mathbb{R}^3$ 
and the total energy in presence of the applied force $f^h\in L^2(S^h, \mathbb{R}^3)$
are given, respectively, by:
\begin{equation*}
\Eh(u^h) = \frac{1}{h}\int_{S^h} W(\nabla u^h), \quad
J^h(u^h) = \Eh(u^h) - \frac{1}{h}\int_{S^h} f^hu^h \qquad \forall u^h\in W^{1,2}(S^h,\mathbb{R}^3).
\end{equation*}
It has been shown \cite{FJMhier} that if $f^h$ scale like $h^\alpha$, then
$\Eh(u^h)$ at approximate minimizers $u^h$ of $J^h$
scale like $h^\beta$, with $\beta = \alpha$ for $0 \le \alpha \le 2$ and $\beta
= 2\alpha -2$ for $\alpha > 2$. The dimension reduction
question in this context consists thus of identifying the $\Gamma$-limits $\mathcal{I}_{\beta, S}$ of 
the rescaled energies sequence $\{\frac{1}{h^{\beta}}\Eh\}_{h\to
  0}$. We stress that, contrary to the curvature-driven shape
formation described in section \ref{sec4}, 
there is no energy quantisation and any scaling exponent $\beta>0$ is viable.

In case of $S\subset\mathbb R^2$ i.e. when $\{S^h\}_{h\to 0}$ is a family of
thin plates, such $\Gamma$-convergence was first established for $\beta = 0$ \cite{LeD-Rao},
and  later \cite{FJMgeo, FJMhier} for all $\beta \geq 2$. 
This last regime corresponds to a rigid  behavior of the elastic material, 
since the limiting deformations are isometries
if $ \beta =2$ (in accordance with the general result in Theorem
\ref{beta2}), or infinitesimal isometries if $\beta>2$ (see, for example, the compactness
analysis in section \ref{sec_compa2}).
One particular case is  $\beta=4$, where the derived limiting theory turns 
out to be  the von K\'arm\'an theory (\ref{vonKar}), then $\beta>4$
with the $\Gamma$-limit as in (\ref{linlin}), and $\beta\in (2,4)$
where the result is effectively included in Theorem \ref{th_Ah}.
We gather these results in Figure \ref{table_plate}, which should be
compared with Figure \ref{table_prestrain} in section \ref{sec_compa2}.

\begin{figure}[htbp]
\centering
\renewcommand{\arraystretch}{1.5}
$\begin{array}{l | l | l | l}
\mbox{\begin{minipage}{1.9cm}   {scaling \\ exponent $\beta$}\end{minipage}} \hspace{-2mm} &
~~~~~\mbox{\begin{minipage}{2.6cm}   {asymptotic \\ expansion of \\  minimizing 
  $u^h_{\mid\omega}$} \end{minipage}} &
~~~~\mbox{\begin{minipage}{2.3cm}  {constraint \\ / regularity}\end{minipage}} &
~~~~~~~~~~~  \Gamma-\mbox{limit }{\mathcal{I}_{\beta, S}} \\ [15pt]\hline
\hspace{-3mm} \begin{array}{l}  {\beta=2} \vspace{-2mm} \\
    \mbox{Kirchhoff} \end{array} & \hspace{-1.5mm}\begin{array}{ll} y(x') \\
  \scriptstyle{\big\{ 3d:~ y(x') + x_3\vec n(x') \big\} } \end{array}
\hspace{-2mm} & \hspace{-2mm} \begin{array}{l} y\in W^{2,2}(\omega,
\mathbb{R}^3)  \\ (\nabla y)^T\nabla y = Id_2 \end{array}
\hspace{-2mm} & c\|(\nabla y)^T\nabla \vec N\|_{\mathcal{Q}_2}^2  \\ [15pt]\hline 
\hspace{-3mm} \begin{array}{l}  {2<\beta<4} \vspace{-2mm} \\
   \mbox {linearised Kirchhoff} \end{array} \hspace{-2mm} & x'+h^{\beta/2-1}v(x')x_3
& \hspace{-2mm} \begin{array}{l} v\in W^{2,2}(\omega,\mathbb{R}) \\ \det\nabla^2v
  =0\end{array} & c\|\nabla^2v\|_{\mathcal{Q}_2}^2 \\  [15pt] \hline 
\hspace{-3mm} \begin{array}{l}   {\beta=4} \vspace{-2mm}  \\
    \mbox{von K\`arm\`an} \end{array} & \hspace{-2mm}\begin{array}{l}x'+hv(x')x_3 \\
  \quad  + h^2w(x')\end{array} &  \hspace{-2mm} \begin{array}{l} v\in
  W^{2,2}(\omega,\mathbb{R}) \\ w\in W^{1,2}(\omega,
  \mathbb{R}^2)\end{array} \hspace{-3mm} & \hspace{-2mm} \begin{array}{l} c_1\|\frac{1}{2}\nabla
  v^{\otimes 2} + (\nabla w)_{sym}\|_{\mathcal{Q}_2}^2 \\ + c_2\|\nabla^2v\|_{\mathcal{Q}_2}^2 \end{array}\\
[15pt] \hline  \hspace{-3mm} \begin{array}{l}   {\beta>4} \vspace{-2mm} \\
   \mbox{linear elasticity} \end{array}  & x'+h^{\beta/2-1}v(x')x_3 & v\in
W^{2,2}(\omega,\mathbb{R}) & c\|\nabla^2v\|_{\mathcal{Q}_2}^2
\end{array}$ 
\caption{The finite hierarchy of $\Gamma$-limits for plates for the
  energy scaling $\beta\geq 2$}
\label{table_plate}
\end{figure}

\subsection{The infinite hierarchy of shell theories and the matching properties}

The first result for the case when $S$ is a surface of arbitrary
geometry was given in \cite{LeD-Rao} as the membrane theory
($\beta=0$) where the limit $\mathcal{I}_{0,S}$ depends only on the stretching 
and shearing produced by the deformation. 
Case $\beta=2$ was analyzed in \cite{FJMM_cr} and proved to reduce
to the flexural shell model, i.e. a geometrically nonlinear 
pure bending, constrained to isometric immersions of $S$.
The energy $\mathcal{I}_{2, S}$ depends then on the change of curvature produced 
by such deformation, in the same spirit of Theorem \ref{beta2}.
 
For $\beta = 4$ the $\Gamma$-limit $\mathcal{I}_{4, S}$,  as shown in \cite{LMP, lemopa2, lemopa_convex},
acts on the first order isometries $V\in \mathcal{V}_1\cap W^{2,2}$ i.e. displacements of $S$
whose covariant derivative is skew-symmetric,
and finite strains $B \in \mbox{cl}_{L^2}\{\mbox{sym}\nabla w; ~ w\in
W^{1,2}(S,\mathbb{R}^3)\}$ (compare the definitions of spaces
$\mathcal{V}_{y_0}$ and $\mathcal{S}_{y_0}$ in section
\ref{sec_compa2}). The limiting energy consists of two terms
corresponding to the stretching (second order change in metric)
and bending (first order change in the second fundamental form
$II=\nabla \vec N$ on $S$) of a family of
deformations $\{\phi^h = \mbox{id} + hV + h^2 w^h\}_{h\to 0}$ of $S$,
which is induced by displacements $V\in\mathcal{V}_1$ and 
$w^h$ satisfying $\lim_{h\to 0}\mathrm{sym} \nabla w^h = B$.
The out-of-plane displacements $v$ present in (\ref{vonKar}) are therefore replaced by the vector fields in
$\mathcal{V}_1$ that are neither normal, nor tangential to $S$,  
but preserving the metric on $S$ up to first order.
For $\beta>4$ the limiting energy
consists  \cite{LMP, lemopa2} only of the bending term and it coincides with the linearly
elastic flexural shell model.

\smallskip

The form of $\mathcal{I}_{\beta, S}$ for all $\beta>2$ and arbitrary
mid-surface $S$ has been conjectured in \cite{lepa2}.
Namely, $\mathcal{I}_{\beta, S}$ acts on the space
of $k$-th order infinitesimal isometries $\mathcal{V}_k$, where $k$ is such that:
$$ \beta\in [\beta_{k+1}, \beta_k) \quad \mbox{ and }  \quad \beta_n= 2+2/n.$$
The space $\mathcal{V}_k$ consists of $k$-tuples $(V_1,\ldots, V_k)$ of displacements 
$V_i:S\rightarrow \mathbb{R}^3$ (with appropriate regularity),
such that the deformations
$\phi^\epsilon = \mathrm{id} + \sum_{i=1}^k \epsilon^i V_i$
preserve the metric on $S$ up to order $\epsilon^k$, i.e.
$(\nabla \phi^\epsilon)^T\nabla \phi^\epsilon - \mbox{Id}_2 = \mathcal{O}(\epsilon^{k+1})$.
Further, setting $\epsilon = h^{\beta/2 - 1}$, we have:
\begin{enumerate}[leftmargin=7mm]
\item[(i)] When $\beta=\beta_{k+1}$ then 
$\mathcal{I}_{\beta, S} = \int_S \mathcal {Q}_2 \left (x, \delta_{k+1}I_S\right)
+ \int_S \mathcal {Q}_2 \left (x, \delta_1 II_S \right)$, 
where $\delta_{k+1}I_S$ is the change of metric
on $S$ of the order $\epsilon^{k+1}$, generated by the family of deformations 
$\{\phi^\epsilon\}_{\epsilon\to 0}$ and $\delta_1 II_S$ is the first order  (i.e. order
$\epsilon$)  change in the second fundamental form $II_S$ of $S$.
\item[(ii)] When $\beta\in (\beta_{k+1}, \beta_k)$ then
$\mathcal{I}_{\beta, S} = \int_S \mathcal {Q}_2 \left (x, \delta_1 II_S \right )$.
\item[(iii)] The constraint of $k$-th order isometry $\mathcal{V}_k$ 
may be relaxed to that of $\mathcal{V}_m$, $m<k$, if $S$ possesses the following 
$m\mapsto k$ matching property. For every $(V_1,\ldots V_m)\in\mathcal{V}_m$ there exist
sequences of corrections $V_{m+1}^\epsilon,\ldots  V_{k}^\epsilon$, uniformly 
bounded in $\epsilon$, such that:
$\tilde \phi^\epsilon = \mathrm{id} + \sum_{i=1}^m \epsilon^i V_i
+ \sum_{i=m+1}^k\epsilon^i V_i^\epsilon$
preserve the metric on $S$ up to order $\epsilon^k$.
\end{enumerate}

\smallskip

The above finding is supported by all the rigorously derived
models. In particular, since plates enjoy the
$2\mapsto\infty$ matching property (i.e., as shown in \cite{FJMhier}, every $W^{1,\infty}\cap
W^{2,2}$ member of $\mathcal{V}_2$ may be matched to an exact
isometry, in the sense of (iii) above), all the plate theories for
$\beta\in(2,4)$ indeed collapse to a single theory (linearized
Kirchhoff model, see Figure \ref{table_plate}).

Further, elliptic (i.e. strictly convex up to
the boundary) surfaces enjoy \cite{lemopa_convex} a matching property of $1\mapsto\infty$, which is 
stronger than for the case of plates.
Namely, on $S$ elliptic and $\mathcal{C}^{3,\alpha}$,
every $V\in\mathcal{V}_1\cap\mathcal{C}^{2,\alpha}(\bar S)$,
possesses a sequence $\{w_\epsilon\}_{\epsilon\to 0}$, equibounded in 
$\mathcal{C}^{2,\alpha}(\bar S, \mathbb{R}^3)$, and such that 
$\phi^\epsilon = \mathrm{id} +\epsilon V + \epsilon^2w_\epsilon$ is an
(exact) isometry for all $\epsilon\ll 1$. 
Regarding the assumed regularity of $V$ (which is higher that
the expected regularity $W^{2,2}$ of a limiting displacement) we note that
the usual mollification techniques do not guarantee the density of smooth 
infinitesimal isometries in $\mathcal{V}_1\cap W^{2,2}$, even for $S\in\mathcal{C}^\infty$.
However, a density result is valid for elliptic $S\in \mathcal{C}^{m+2,\alpha} $, that is: 
for every $V\in\mathcal{V}_1\cap W^{2,2}$ there exists a sequence 
$\{V_n\in\mathcal{V}_1\cap\mathcal{C}^{m,\alpha}(\bar S,\mathbb{R}^3)\}_{n\to\infty}$ 
such that: $\lim_{n\to\infty} \|V_n - V\|_{W^{2,2}(S)} = 0.$
The proof of the quoted results adapts techniques used for 
showing immersability of all positive curvature metrics on a sphere \cite{HH}.
As a consequence, for elliptic surfaces with sufficient regularity the $\Gamma$-limit of 
the nonlinear elastic energies $h^{-\beta}\Eh$ for any scaling regime
$\beta>2$ is given by the bending functional constrained to the first
order isometries, as in the case $\beta>4$.

In \cite{holepa, horn} a matching and density
properties of isometries on developable surfaces without affine
regions, has been proved. Namely, on such $S$ of regularity $\mathcal{C}^{2k,
  1}$, every $V\in\mathcal{V}_1\cap \mathcal{C}^{2k-1,1}$ enjoys
$1\mapsto k$ matching property. Further, the space $\mathcal{V}_1\cap
\mathcal{C}^{2k-1}$ is dense in $\mathcal{V}_1$. 
The implication for elasticity of thin shells with smooth developable mid-surface
is that, again, the only small slope theory is the linear theory; a
developable shell transitions directly from the linear regime to fully
nonlinear bending if the applied forces are adequately increased.
While the von K\'arm\'an theory describes buckling of thin plates, the equivalent
variationally correct theory for developable shells is the purely nonlinear bending.
It is worth to notice that the class of developable shells includes smooth
cylinders which are ubiquitous in nature and technology over a range of length scales. An
example of a recently discovered structure is carbon nanotubes, i.e. molecular-scale tubes of
graphitic carbon with outstanding rigidity properties: they are among the stiffest materials
in terms of the tensile strength and elastic modulus, but they easily buckle under compressive,
torsional or bending stress.

\begin{problem}
Investigate the matching properties for other types of surfaces.
\end{problem}

\section{Future directions}\label{sec_fut}

Our review on the mathematical aspects of the morphogenesis and pattern formation in thin sheets has focused on  low-dimensional shapes that arise from inhomogeneous growth and/or boundary conditions and constraints.   From a biological perspective, understanding how growth leads to shape is only half the problem. A true understanding of morphogenesis requires to also understand how shape feeds back to growth, to ultimately regulate shape and thus enable function. From a technological perspective, an equally interesting problem is the inverse problem: how should one prescribe the growth patterns in order to be able to convert a flat sheet into a complex landscape, a flower or even a face?

From an artistic perspective, a natural generalization of the questions on the smoothness of and in pattern-forming elastic surfaces, is that posed by the ancient Sino-Japanese paper arts of origami and kirigami (from the japanese: Oru=fold, Kiri=cut, Kami=paper): what are the limits to the shapes that one can construct with sharp folds and cuts, that violate smoothness along cuts and creases (either straight or curved)? Artists have long known how to fold a sheet into a crane, a man or a dragon, and how to use cuts to articulate a sheet so that it can be made into a pop-up castle or a rose.  How can one quantify these art forms as inverse problems in discrete geometry and topology? We touch on each of these three problems briefly to highlight recent progress, and the many open problems that remain.

\subsection{Developmental feedback from shape to growth} In a biological context, there is increasing evidence for a mechanical feedback loop linking shape back to growth \cite{Shraiman,Shraiman2}, i.e. the growth tensors associated with causing shape are themselves affected by shape. To quantify how growth patterns change in response to shape in space and time with (unknown) kernels that characterize the nature of this feedback, one must turn to experiments. Nevertheless, it might still be useful to study simple feedback laws to understand their mathematical consequences as has been recently attempted in the context of controlling the bacterial shapes \cite{Santangelo}. A minimal example of a local model, incorporating mechanical feedback in tissue growth (in such instances as leaves and epithelial tissues), that closes the equations (\ref{Goveq0}), takes the form:
\begin{equation}\label{Goveq5}
\begin{split}
& \alpha_v\Delta ({\rm tr}\,\dot{\mathbf{s}}) = -\alpha\Delta ({\rm tr}\, \mathbf{s}) -\frac{\alpha}{2} \det (\boldsymbol \kappa_0- \boldsymbol\kappa)-  \Delta [{\rm tr}\,(\boldsymbol \sigma_0-\boldsymbol \sigma)]\\
& -\beta_v \Delta ({\rm tr}\, \dot{\mathbf{b}})=-\beta \Delta ({\rm tr}\, \mathbf{b}) -\beta\Delta [{\rm tr} (\boldsymbol \kappa_0- \boldsymbol\kappa)]+ {\rm tr}\, [(\boldsymbol \sigma_0- \boldsymbol\sigma) (\boldsymbol \kappa_0-\boldsymbol \kappa)],
\end{split}
\end{equation}
for the dynamics of the in-plane growth and curvature tensors $\mathbf{s}, \mathbf{b}$, respectively. Here, the terms $\boldsymbol \sigma_0,\boldsymbol\kappa_0$ denote the threshold homeostatic values of the stress tensor $\boldsymbol\sigma$ and the curvature tensor $\boldsymbol\kappa$ that the tissue aims to achieve, and the various prefactors are as defined in the introductory section, except for $\alpha_v, \beta_v$ which are the stretching viscosity and bending viscosity respectively, with $\alpha_v/\alpha = \tau_S, \beta_v/\beta = \tau_B$ being the time scale for the relaxation of in-plane and out-of-plane growth. We note that these equations are linear in $\mathbf{s},\mathbf{b}$ and thus likely to be valid only in the neighborhood of the homeostatic stress and curvature. 

\begin{problem}
The system (\ref{Goveq5}) is geometrically nonlinear. What are the conditions for its dynamic stability and control of the equilibrium states, that result from inhomogeneous and anisotropic growth?
\end{problem}

Other possible descriptions were suggested in \cite{jed1, jed2, jed3, bre_lew, lemucha}. In particular, the paper \cite{bre_lew} has recently introduced a free boundary problem for a system of pdes modeling growth. There, a morphogen controlling volume growth and produced by specific
cells, was assumed to be diffused and absorbed throughout the domain, whose geometric shape was in turn determined by the instantaneous minimization of an elastic deformation energy, subject to a constraint on the volumetric growth. For an initial domain with $\mathcal{C}^{2,\alpha}$ regular
boundary, it has been establishes the local existence and uniqueness of a classical
solution, up to a rigid motion.



\subsection{Inverse problems in morphogenesis}
\begin{figure}
   \centering
    \includegraphics[width=\textwidth]{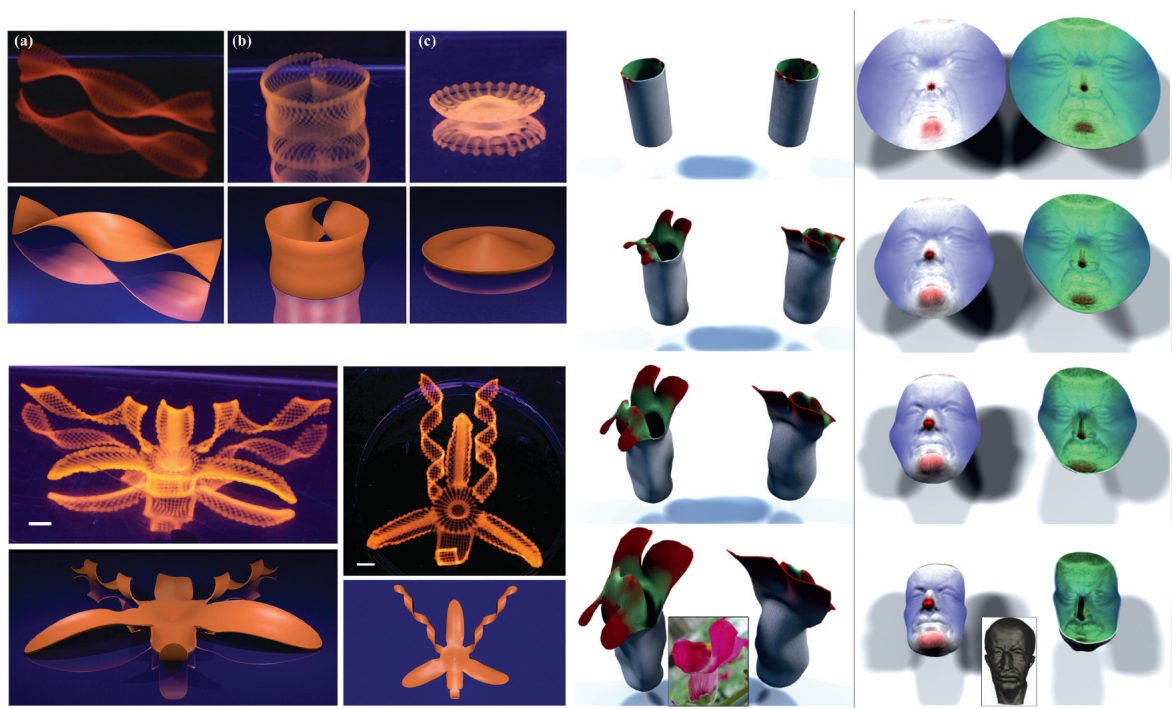}
   \vskip -.6in
   \caption{Examples of solutions of inverse problems in morphogenesis - to program the metric structure and thence create complex shapes from flat sheets. On the left \cite{22} are shown experiments with 3d printed gel structures that swell in a good solvent, along with representative numerical solutions that are based on minimizing the energy (\ref{Ih}). On the right \cite[Copyright (2021) National Academy of Sciences]{RVM} we see results of the solution of inverse problems to grow a flower from a bilayer cylindrical shell, and a face from a circular bilayer disk.}
  \label{fig7}
\end{figure}

With the advent of additive manufacturing methods such as 3d and 4d printing (to account for variations in space and time), it has now become possible to print planar patterns of responsive inks that swell or shrink inhomogeneously when subject to light, pH, humidity etc. thus causing them to bend and twist out of the plane \cite{22,boley2019shape}. Understanding how to design the ink materials and the geometric print paths to vary the density and anisotropy of the print patterns in a monolayer or a bilayer is critical to enable functional patterns. This inverse design problem requires the specification of the first and second fundamental form which will not generally be compatible with a strain-free final shape. Recent work in this area \cite{RVM} shows that a way around this is to use a bilayer with independent control over the two layers and leads to results such as those shown in Figure \ref{fig7}. A related class of design problems in solid mechanics, leading to a variation on the classical question of equi-dimensional embeddability of Riemannian manifolds has been addressed in \cite{ALP}.

\subsection{Discrete problems: origami and kirigami}

\begin{figure}
   \centering
    \includegraphics[width=\textwidth]{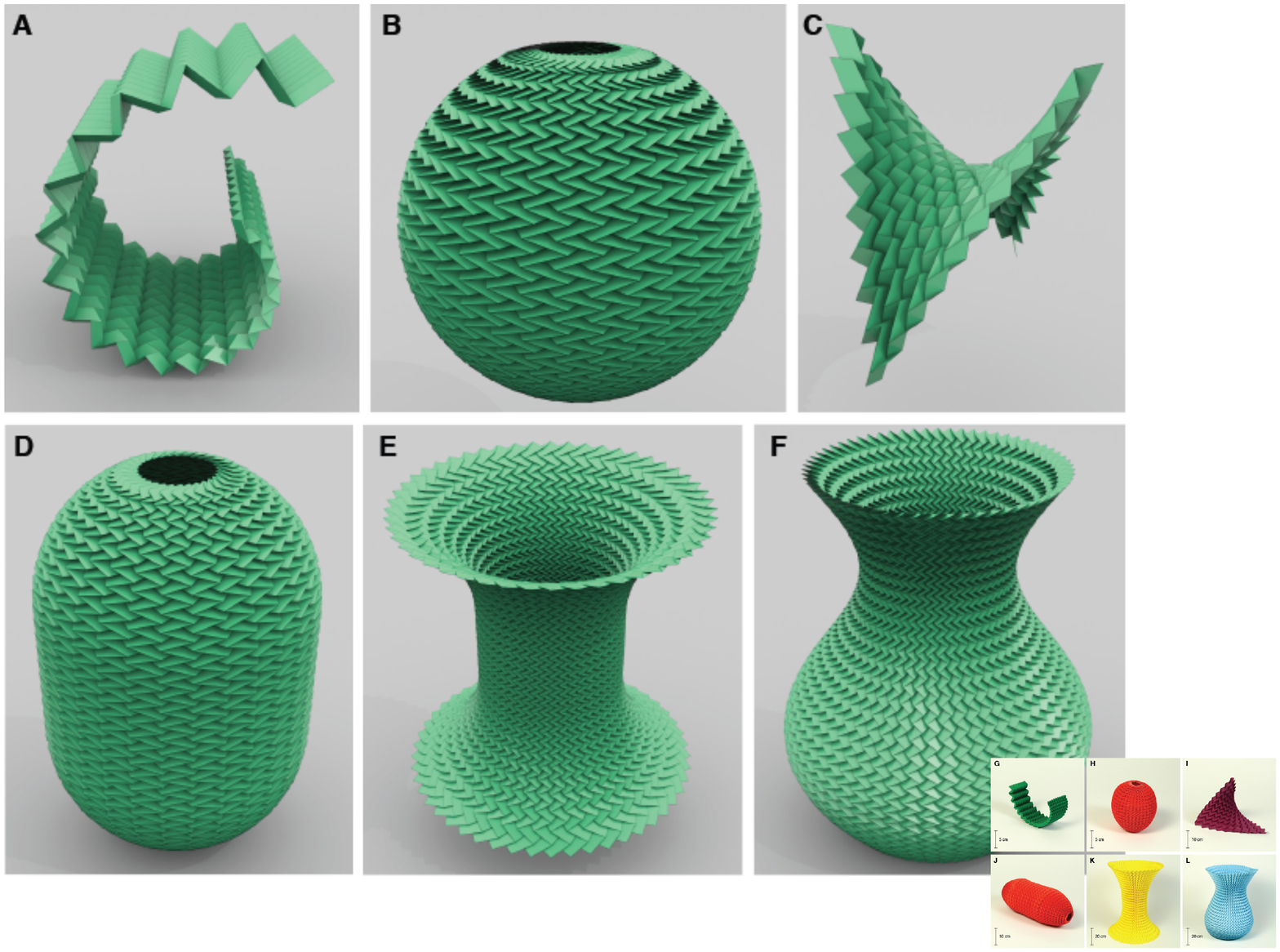}
   \vskip -.2in
   \caption{While this paper has focused on surfaces with varying degree of smoothness, an interesting 
new avenue for exploration is that of discrete surfaces \cite{Dudte_origami} that have strongly creased 
regions, seen for example in Origami.  
   }
  \label{fig8}
\end{figure}

Origami is the art of folding paper along sharp creases to create complex three-dimensional shapes, and thus more amenable to the methods of discrete geometry. A natural question here is that of designing the number, location and orientation of folds on a flat sheet of paper and prescribing the order of folding to achieve a given target shape. For a prescribed fold topology, e.g. that of 4-coordinated vertices, geometric rules that quantify the constraints of local length, angle and area preservation allow one pose the inverse problem of fold design as a constrained optimization problem\cite{Dudte_origami,D3, D2}. Then, given reasonable initial states, one can determine the folding patterns to achieve target shapes (see Figure \ref{fig8}) that are realized as spatially modulated patterns of a simple periodic and uniform tiling yield approximations to given surfaces of constant or varying curvature, and corroborated using experiments with paper. The difficulty of realizing these geometric structures may be assessed by quantifying the energetic barrier that separates the metastable flat and folded states. The trade-off between the accuracy to which the pattern conforms to the target surface and the effort associated with creating finer folds, can also be characterized \cite{DO}. However, there are a host of mathematical problems that remain open. These include the presence (or absence) of impossibility theorems on what shapes can or cannot be achieved using folds in a sheet of paper, and the consequences of fold topology on the resulting shapes. 

\begin{problem}
How can one control the rigidity of a randomly origamized sheet, as the number of random creases is gradually increased, and subject to the geometric rules that the creases must satisfy at every vertex, i.e. the sum of all angles must add up to $2\pi^c$, and that alternate angles must add up to $\pi^c$ \cite{Dudte_origami}?
\end{problem}

Kirigami is the art of cutting paper to make it articulated and deployable.
The mechanical response of a kirigami
sheet when it is pulled at its ends is enabled and limited by the
presence of cuts that serve to guide the possible
non-planar deformations. Recently, this ability has become the
inspiration for a new class of mechanical metamaterials
\cite{Kirigami-review, Metamaterial-review}. The
geometrical and topological properties of the slender sheet-like
structures, irrespective of their material constituents, were exploited to discuss
functional structures on scales ranging from the nanometric
\cite{Mceuen} to centimetric and beyond \cite{Bertoldi, M1, M2}.  

A combination of physical and
numerical experiments can be used to characterize the geometric
mechanics of kirigamized sheets as a function of the number, size, and
orientation of cuts \cite{Kirigami-jointpaper}.  In particular, of interest is understanding how the varying
of the the shortest path
between points at which forces are applied, influences the shape of the deployment of the trajectory of a sheet as well as how to control its compliance across orders of magnitude.

Mathematically, these questions are related to the nature and form of
geodesics in the Euclidean plane with linear obstructions (cuts), and
the nature and form of isometric immersions of the sheet with cuts
when it can be folded on itself.  In \cite{Kiri-HLM}, a constructive proof
has been provided that the geodesic connecting any two points in the sheet is piecewise
polygonal, and that the family of all such geodesics can
be simultaneously rectified into a straight line by flat-folding the sheet so that its
configuration is a (non-unique) piecewise affine isometric immersion.  

\begin{problem}
Study the structure of geodesics in the kirigamized sheet as the number of random cuts increases to infinity, and under various assumptions on the cuts distribution. What is the Hausdorff dimension of the limiting paths?
\end{problem}

\end{document}